\documentstyle[12pt]{article}
\begin{document}
\title{Differentiable functions of quaternion variables.}
\author{S.V. L\"udkovsky, F. van Oystaeyen.}
\date{16.08.2002}
\maketitle
\begin{abstract}
We investigate differentiability of functions defined on regions of the
real quaternion field and obtain a noncommutative version
of the Cauchy-Riemann conditions. Then we study
the noncommutative analog of the Cauchy integral
as well as criteria for functions of a quternion variable
to be analytic. In particular, the quaternionic
exponential and logarithmic functions are being considered.
Main results include quaternion versions of Hurwitz' theorem,
Mittag-Leffler's theorem and Weierstrass' theorem.
\end{abstract}
\section{Introduction}
The noncommutativity of the quaternion field $\bf H$
obstructs the immediate application of the theory
of analytic and meromoprhic complex functions to functions of
quaternion arguments.
The latter may be thought of as functions of two noncommuting complex
variables, but we shall adopt matrix notation representing
the standard generators of the quaternions by their Pauli-matrices.
This allows a rather elegant introduction
of differentiable functions of a quaternion variable,
integrals of functions along curves in $\bf H$, residues of
a function,... .The new results contained in this paper
provide noncommutative analogs of the Cauchy-Riemann
conditions for superdifferentiable functions as well
as basic properties of the related noncommutative
integrals, the argument principle, etc. ...
The quaternionic residue theory depends on the definition
and description of the exponential and logarithm functions
of quaternion variables. An explicit description
of the exponential is obtained in Proposition 3.2 allowing
to view it as an epimorphism from a set of imaginary quaternions
to the three dimensional
quaternionic unit sphere. The relation between the quaternion
version of holomorphicity and local analyticity is investigated,
in particular we obtain in Theorem 3.10 that for a continuous function
on an open subset $U$ of $\bf H$, the property of being locally
analytic follows from the integral holomorphicity.
Section 3 also contains the quaternionic version of the classical
theorems of Cauchy, Liouville and Morera.
Although the analytic theory of functions of a quaternion
variable, or more general of functions of noncommuting variables
with suitable "commutation" rules, has an interest in its own right,
we were more motivated by the connection with noncommutative geometry,
the analytic structure induced on $\bf H$ modulo a $\bf Z$-lattice,
and the quaternionic version of arithmetical functions like
the zeta-function.  We hope to return to these applications in
forthcoming work.
\par Though some results in noncommutative geometry are concerned with
function families \cite{berez,connes,dewitt,khren,oystaey}
they are rather general and do not take into account
the particular quaternion case and its specific features.
It is necessary to note, that we use a weaker superdifferentiability
condition, compared to, for example, \cite{berez,dewitt,khren}.
Traditionally one uses the condition, that a
right derivative is right superlinear on a superalgebra, which
causes severe restrictions on these classes of functions
(see Theorems I.1.4 and I.2.3 in \cite{khren}).
This is too restrictive in the particular quaternion case
as it does not permit to describe an $\bf H$-algebra of quaternion
holomorphic functions on an open subset $U$ in $\bf H^n$
extending that of complex holomorphic functions.
We have withdrawn the condition of right superlinearity
of a superderivative on a superalgebra, supposing only that it is
additive on $\bf H^n$ and $\bf R$-homogeneous.
Nevertheless, it also satisfies distributivity and
associativity laws relative to the multiplication from the right
on (scalar) quaternions $\lambda \in \bf H$
and there are also distributivity and associativity
laws relative to a left multiplication on $\lambda \in \bf H$
(see \S 2.1).
That is, we have considered Frech\'et differentiable
functions on the Euclidean space $\bf R^{4n}$
with some additional conditions on increments of
functions, taking into account a superalgebra structure.
This permits to encompass classes of all
analytic functions on a region in $\bf H^n$, in particular,
all polynomial functions.  Moreover, this approach
permits to investigate an analog of functions having
Laurent series expansions. We have proved, that for each complex
holomorphic function $f$ on a region $V$ open in $\bf C^n$
there exists a quaternion holomorphic function $F$ on a
suitable region $U$ in $\bf H^n$ such that a restriction of $F$ on $V$ 
coincides with $f$.
The theory of complex holomorphic functions turns out to be
rather different from a theory of quaternion
holomorphic functions. The quaternion field $\bf H$ has
nontrivial algebraic structure and identities, so there are different ways
to define not only function spaces, but also their differentiations.
A differentiation is not only analytic, it also
has algebraic properties. In some sense the notion of the family
of all quaternion holomorphic
functions unifies together complex holomorphic and antiholomorphic
functions.
On the other hand, weaker differentiability, for example,
"pointwise" as defined in the classical case
by Gat$\hat o$ \cite{kolfom} yields a too poor
algebraic structure of function spaces not taking
into account the gradation of a superalgebra.
\par In previous works \cite{lulgcm,lulsqm,lustptg,lupm}
the first author investigated loop and diffeomorphism groups
of complex manifolfds and quasi-invariant measures 
and stochastic processes on them. Complex manifolds also
have the structure of supermanifolds, since the field
$\bf C$ can be considered as a graded algebra over $\bf R$.
The graded structure of the quaternion field over the reals
is more complicated. Conceivably, the investigations
in this work allow to continue
this work for quaternion manifolds, loop
and diffeomorphism groups of these manifolds, quasi-invariant measures
and stochastic processes on them,
as well as their associated representations
including irreducible ones.
\section{Differentiability of functions of quaternion variables}
\par To avoid misunderstandings we first introduce notations.
We write $\bf H$ for the skewfield of quaternions over
the real field $\bf R$. This skewfield can be represented
as a subring of the ring ${\bf M}_2({\bf C})$ of all
$2\times 2$ complex matrices by representing the classical quaternion
basis $1$, $i$, $j$, $k$, by the Pauli-matrices
$I$, $J$, $K$, $L$ defined as follows:
$$I={{1\quad 0}\choose {0\quad 1}};\quad J={{i\quad 0}\choose {0\quad -i}};
\quad K={{0\quad 1}\choose {-1\quad 0}};
\quad L={{0\quad i}\choose {i\quad 0}},$$
where $i=(-1)^{1/2}$.
Hence, each quaternion $z$ is written
as a $2\times 2$ matrix over $\bf C$ having matrix elements
$z_{1,1}={\bar z}_{2,2}=:t$, $z_{1,2}=-{\bar z}_{2,1}=:u$,
where $t$ and $u\in \bf C$ such that $t=v+iw$ and $u=x+iy$,
$v$, $w$, $x$ and $y$ are in the field $\bf R$ of real numbers.
The quaternion skewfield $\bf H$ has an anti-automorphism $\eta $
of order two induced in $\bf H$ by the Hermite conjugation
in ${\bf M}_2({\bf C})$, that is,
$\eta: z\mapsto {\tilde z}$, where ${\tilde z}_{1,1}=\bar t$ and
${\tilde z}_{1,2}=-u$. There is a norm in $\bf H$
such that $|z|=(|t|^2+|u|^2)^{1/2}$, hence $det(z)=|z|^2$
and ${\tilde z}=|z|^2z^{-1}$. The noncommutative field $\bf H$
is the $\bf Z_2$-graded $\bf R$-algebra ${\bf H}={\bf H}_0+{\bf H}_1$, where
elements of ${\bf H}_0$ are $\underline {even}$ 
and elements of ${\bf H}_1$ are $\underline {odd}$
(see, for example, \cite{cohn,kosshaf,waerd}).
\par In view of noncommutativity of $\bf H$ a polynomial function
$P: U\to \bf H$ may have several different representations
$$\mbox{\v{P}}(z,\tilde z)=\sum_kb_{k,1}{\hat z}^{k_1}...b_{k,m}{\hat z}^{k_m},$$
where $b_{k,j}\in \bf H$ are constants, $k=(k_1,...,k_m)$,
$m\in \bf N$, $k_j=(k_{j,1},...,k_{j,2n})$, $k_{j,l}\in \bf Z$,
${\hat z}^{k_j}:=\mbox{ }^1z^{k_{j,1}}\mbox{ }^1{\tilde z}^{k_{j,n+1}}...
\mbox{ }^nz^{k_{j,n}}\mbox{ }^n{\tilde z}^{k_{j,2n}},$
$\mbox{ }^lz^0:=1$, $\mbox{ }^l{\tilde z}^0=1$, $U$ is an open subset
of $\bf H^n$. Each term  $b_{k,1}{\hat z}^{k_1}...b_{k,m}{\hat z}^{k_m}=:
\omega (b_k,z,\tilde z)\ne 0$ we consider as a word of length
$\xi (\omega )=\sum_{j,l}\delta (k_{j,l}) +\sum_j\kappa (b_{k,j})$,
where $\delta (k_{j,l})=0$ for $k_{j,l}=0$ and $\delta (k_{j,l})=1$
for $k_{j,l}\ne 0$, $\kappa (b_{k,j})=j$ for $b_{k,j}=1$,
$\kappa (b_{k,j})=j+1$ for $b_{k,j}\in {\bf H}\setminus \{ 0,1 \} $.
A polynomial $P$ is considered as a phrase $\mbox{\v{P}}$ of a length
$\xi (\mbox{\v{P}}):=\sum_k\xi (\omega (b_k,z,\tilde z))$.
Using multiplication of constants in $\bf H$, commutativity
of $vI$ with each $\mbox{ }^lz$ and $\mbox{ }^l\tilde z$,
and $\mbox{ }^lz^a\mbox{ }^lz^b=\mbox{ }^lz^{a+b}$
and $\mbox{ }^l{\tilde z}^a\mbox{ }^l{\tilde z}^b=
\mbox{ }^l{\tilde z}^{a+b}$, $\mbox{ }^lz\mbox{ }^l{\tilde z}=
\mbox{ }^l{\tilde z}\mbox{ }^lz$, it is possible to consider
representations of $P$ as phrases $\mbox{\v{P}}$ of a minimal lenght
$\xi (\mbox{\v{P}})$. We choose one such $\mbox{\v{P}}$ of a minimal lenght.
If $f: U\to \bf H$ is a function presented by a convergent by
$z$ and $\tilde z$ series $f(z,\tilde z)=\sum_nP_n(z,\tilde z)$,
where $P_n(vz,v\tilde z)=v^nP_n(z,\tilde z)$ for each $v\in \bf R$
is a $\bf R$-homogeneous polynomial, $n\in \bf Z$,
then we consider among all representations of $f$
such for which $\xi (\mbox{\v{P}}_n)$ is minimal for each
$n\in \bf Z$. We may use this convention separately
for families of functions $f$ having
$(a)$ $z$-series decompositions, $(b)$ $\tilde z$-series decompositions,
$(c)$ $(z,\tilde z)$-series decompositions (that is,
by indicated variables).
The corresponding families of locally analytic functions on $U$
are denoted by $C^{\omega }_z(U,{\bf H})$,
$C^{\omega }_{\tilde z}(U,{\bf H})$,
$C^{\omega }_{z,\tilde z}(U,{\bf H}).$  If each $P_n$ for $f$
has a decomposition of a particular left type
$$\mbox{\v{P}}(z,\tilde z)=\sum_{k,p}b_{k,p}z^k{\tilde z}^p,$$
where $0\le k, p\in \bf Z$, $b_{k,p}\in \bf H$, then
the space of all such locally analytic functions on $U$
is denoted by $\mbox{ }_lC^{\omega }_{z,\tilde z}(U,{\bf H}),$
for $z$-series or $\tilde z$-series decompositions only
the corresponding spaces are denoted by
$\mbox{ }_lC^{\omega }_z(U,{\bf H})$  and
$\mbox{ }_lC^{\omega }_{\tilde z}(U,{\bf H})$  respectively.
They are proper subspaces of that of given above.
Spaces of locally analytic functions $f$ having right type
decompositions for each $P_n$
$$\mbox{\v{P}}(z,\tilde z)=\sum_{k,p}z^k{\tilde z}^pb_{k,p}$$
are denoted by
$\mbox{ }_rC^{\omega }_{z,\tilde z}(U,{\bf H}),$ etc.
\par {\bf 2.1. Definition.} Consider
an open region $U$ in $\bf H^n$, the $n$-fold product of copies
of $\bf H$, and let $f: U\to \bf H$
be a function. Then $f$ is said to be (right) superdifferentiable
at a point $(\mbox{ }^1z,...,\mbox{ }^nz)=e_1\mbox{ }^1z+...
+e_n\mbox{ }^nz\in  U$ (with respect to a chosen (right) $\bf H$-basis
for $\bf H^n$, $ \{ e_1,...,e_n \} $),
if it can be written in the form
$$f(z+h)=f(z)+\sum_{j=1}^nA_j\mbox{ }^jh
+\epsilon (h)|h|$$
for each $h\in \bf H^n$ such that $z+h\in U$, where
$A_j\in \bf H$ for each $j=1,...,n$ and $A_j$ is denoted
by $\partial f(z)/\partial \mbox{ }^jz$, that is,
there exists a (right) derivative $f'(z)$ such that
a (right) differential is given by
$$D_zf(z).h:=f'(z).h:=
\sum_{j=1}^n (\partial f(z)/\partial \mbox{ }^jz)\mbox{ }^jh,$$
where $\epsilon (h)$ is a function continuous at zero such that
$\epsilon (0)=0$, \\
$e_j=(0,...,0,1,0,...,0)$ is the vector
in $\bf H^n$ with $1$ on $j$-th place,
$$D_zf(z).h=:(Df)(z;h)$$
such that $(Df)(z;h)$ is additive in
$h$ and $\bf R$-homogeneous, that is,
$$(Df)(z;h_1+h_2)=(Df)(h_1)+(Df)(h_2)\mbox{ and }
(Df)(z;vh)=v(Df)(z;h)$$
for each $h_1$, $h_2$ and $h\in \bf H^n$, $v\in \bf R$.
There are imposed conditions:
$$D_{\tilde z}z=0,\quad D_z{\tilde z}=0,\quad (D_zz).h=h,
\quad D_z1=0,\quad (D_{\tilde z}{\tilde z}).h=h,
\quad D_{\tilde z}1=0,$$
$$\mbox{also } (D_z(fg)).h=((D_zf).h)g+f(D_zg).h$$
for a product of two supedifferentiable functions $f$ and $g$
and each $h\in \bf H^n$. We also have distributivity and
associativity laws relative to multiplication from the right
by (scalar) quaternions $\lambda \in \bf H$:
$$(D(f+g))(z;h\lambda )=(Df)(z;h\lambda )+(Dg)(z;h\lambda ),$$
$$(Df)(z;h(\lambda _1+\lambda _2))=(Df)(z;h\lambda _1)+
(Df)(z;h\lambda _2), $$
$$(Df)(z;(h\lambda _1)\lambda _2)=(Df)(z;h(\lambda _1\lambda _2))$$
for each superdifferentiable functions
$f$ and $g$ at $z$ and each $\lambda $, $\lambda _1$ and
$\lambda _2\in \bf H$.  There are also left distributive
and associative laws:
$$(D\lambda (f+g))(z;h)=\lambda (Df)(z;h)+\lambda (Dg)(z;h),$$
$$(D(\lambda _1+\lambda _2)f)(z;h)=\lambda _1(Df)(z;h)+
\lambda _2(Df)(z;h),$$
$$(D(\lambda _1\lambda _2)f)(z;h)=\lambda _1(D\lambda _2f)(z;h).$$
\par That is, we consider Frech\'et differentiable
functions on the Euclidean space $\bf R^{4n}$
with some additional conditions on increments of
functions, taking into account a superalgebra structure.
Quite analogously
we defined the notion of (right) superdifferentiability by
$\tilde z$ and by their pair $(z,\tilde z).$
\par {\bf Notation.} We write $f$ as a $2\times 2$ complex matrix
with entries $f_{i,j}$ such that $f_{i,j}=g_{i,j}+ih_{i,j}$
and $g_{i,j}$, $h_{i,j}$ being real-valued functions.
For $n=1$ we also write $\mbox{ }^1z$ without its superscript. 
We may write a function $f(z,\tilde z)$ in variables
$(v,w,x,y)$ as $F(v,w,x,y)=f\circ \sigma (v,w,x,y)$, where
$\sigma (\mbox{ }^lv,\mbox{ }^lw,\mbox{ }^lx,\mbox{ }^ly)=
(\mbox{ }^lz,\mbox{ }^l{\tilde z})$ is a bijective mapping.
\par {\bf 2.2. Proposition.} {\it A function $f: U\to \bf H$ is
(right) superdifferentiable at a point $a\in U$ if and only if
$F$ is Frech\'et differentiable at $a$ and
$$(2.1)\quad D_{\tilde z}f(z)|_{z=a}=0.$$
If in addition $f'(a)$ is right superlinear on the superalgebra
$\bf H^n$, then $f$ is superdifferentiable at $a$ if and only if
$F$ is Frech\'et differentiable at $a$ and satisfies
the following equations:
$$(2.2)\quad
\partial G_{1,1}/\partial \mbox{ }^jv=\partial H_{1,1}/\partial \mbox{ }^jw,
\quad
\partial G_{1,1}/\partial \mbox{ }^jw=-\partial H_{1,1}/\partial
\mbox{ }^jv,$$
$$\partial G_{1,2}/\partial \mbox{ }^jv=-\partial H_{1,2}/\partial
\mbox{ }^jw,  \quad
\partial G_{1,2}/\partial \mbox{ }^jw=\partial H_{1,2}/\partial
\mbox{ }^jv,$$
$$\partial G_{1,1}/\partial \mbox{ }^jw=-\partial H_{1,2}/
\partial \mbox{ }^jx,
\quad
\partial G_{1,1}/\partial \mbox{ }^jx=\partial H_{1,2}/
\partial \mbox{ }^jw,$$
$$\partial G_{1,2}/\partial \mbox{ }^jw=-\partial H_{1,1}/
\partial \mbox{ }^jx,
\quad
\partial G_{1,2}/\partial \mbox{ }^jx=\partial H_{1,1}/
\partial \mbox{ }^jw,$$
$$\partial G_{1,1}/\partial \mbox{ }^jx=-\partial H_{1,1}/
\partial \mbox{ }^jy,
\quad
\partial G_{1,1}/\partial \mbox{ }^jy=\partial H_{1,1}/
\partial \mbox{ }^jx,$$
$$\partial G_{1,2}/\partial \mbox{ }^jx=\partial H_{1,2}/
\partial \mbox{ }^jy,
\quad
\partial G_{1,2}/\partial \mbox{ }^jy=-
\partial H_{1,2}/\partial \mbox{ }^jx,$$
or shortly in matrix notation:
$$(2.3)\quad \partial F/\partial \mbox{ }^jv=
(\partial F/\partial \mbox{ }^jw)J^{-1}=
(\partial F/\partial \mbox{ }^jx)K^{-1}=
(\partial F/\partial \mbox{ }^jy)L^{-1}$$
for each $j=1,...,n$.}
\par {\bf Proof.}
Verify that $z$ and $\tilde z$ are independent variables.
Suppose contrary that there exists $\gamma \in \bf H$ such that
$z+\gamma \tilde z=0$ for each $z\in \bf H$.
This is equivalent to a system of two linear equations
$\gamma _{1,1}{\bar t}-\gamma _{1,2}u=-t$ and
$\gamma _{1,1}{\bar u}+\gamma _{1,2}t=-u$.
If $z\ne 0$, then $\gamma _{1,1}=-(t^2+u^2)/(|t|^2+|u|^2)$ and
$\gamma _{1,2}=(t{\bar u}-{\bar t}u)/(|t|^2+|u|^2)$, hence
$\partial \gamma /\partial t\ne 0$ and $\partial \gamma /\partial u\ne 0$.
Therefore, there is not any $\gamma \in \bf H$ such that
$z+\gamma \tilde z=0$ for each $z\in \bf H$.
\par For each canonical closed compact set $U$ in $\bf H$
the set of all polynomial by $z$ and $\tilde z$
functions is dense in the space of all continuous on $U$
Frech\'et differentiable functions on $Int (U)$.
In particular functions of the form of series $f=\sum \mbox{ }_{l_1}f...
\mbox{ }_{l_n}f$ converging on $U$ together with its
superdifferential on $Int (U)$
such that each $\mbox{ }_lf$ is (right) superlinearly
superdifferentiable on $Int (U)$ relative to the superalgebra
$\bf H^n$ is dense in the $\bf R$-linear space of (right) superdifferentiable
functions.
From conditions of \S 2.1 it follows, that the superdifferentiability
conditions are defined uniquely on space of polynomials.
Such that the superdifferentiability of a polynomial $P$ on $U$
means that it is expressible through a sum of products of $\mbox{ }^jz$
and constants from $\bf H$, that is, without terms containing $\tilde z$.
Suppose that $f$ is superdifferentiable at a point $a$.
To each $f'(z)$ there corresponds a $\bf R$-linear
operator on the Euclidean space $\bf R^{4n}$. Moreover,
we have the distributivity and associativity
laws for $(Df)(z;h)$ relative to the right multiplication
on quaternions $\lambda \in \bf H$ (see \S 2.1).
Then $f(a+h)-f(a)=D_af(a,\tilde a).h+D_{\tilde a}f(a,\tilde a).{\tilde h}+
\epsilon (h)|h|$ $=D_af(a,\tilde a).h +\epsilon (h)|h|,$
where $\epsilon (h)$ is continuous by $h$ and $\epsilon (0)=0$,
therefore, $D_{\tilde a}f=0$.
Vice versa, if $F$ is Frech\'et differentiable, then
expressing $vI$, $wJ$, $xK$ and $yL$ through linear combinations
of $z$ and $\tilde z$ with constant coefficients we get
the increment of $f$ as above which is independent of $\tilde h$
if and only if $D_{\tilde a}f=0$.
\par Consider now the particular case, when
$f'$ is right superlinear on the superalgebra $\bf H^n$.
In this case $f'(a)$ is right $\bf H$-linear.
Using the definition of the right
superderivative and that there is a bijective correspondence
between $z$ and $(v,w,x,y)$ we consider a function
$f=f(z,\tilde z)=F(v,w,x,y)$ (right) superdifferentiable by
$z$ and $\tilde z$, hence it is also differentiable by
$(v,w,x,y)=(b_1,...,b_4)$ and we obtain the expressions:
$$\partial F/\partial \mbox{ }^jb_l=(\partial F/\partial \mbox{ }^jz).
\partial \mbox{ }^jz/\partial \mbox{ }^jb_l+
(\partial F/\partial \mbox{ }^j{\tilde z}).
\partial \mbox{ }^j{\tilde z}/\partial \mbox{ }^jb_l,$$
since $\partial \mbox{ }^jz/\partial \mbox{ }^kb_l=0$
and $\partial \mbox{ }^j{\tilde z}/\partial \mbox{ }^kb_l=0$
for each $k\ne l$.
From $D_{\tilde z}f=0$ and
$\partial \mbox{ }^jz/\partial \mbox{ }^jv=I,$
$\partial \mbox{ }^jz/\partial \mbox{ }^jw=J,$
$\partial \mbox{ }^jz/\partial \mbox{ }^jx=K,$
$\partial \mbox{ }^jz/\partial \mbox{ }^jy=L$
we get Equations $(2.3)$.
Substituting $J^{-1}$ for the equation with pair of variables
$(v,w)$ in $(2.3)$ we get
$\partial F_{1,1}/\partial \mbox{ }^jv=
-i\partial F_{1,1}/\partial \mbox{ }^jw$ and
$\partial F_{1,2}/\partial \mbox{ }^jv=
i\partial F_{1,2}/\partial \mbox{ }^jw$,
substituting $K^{-1}J=L$ for Equation $(2.3)$ with
the pair of variables $(w,x)$ we get:
$\partial F_{1,1}/\partial \mbox{ }^jw=
i\partial F_{1,2}/\partial \mbox{ }^jx$ and
$\partial F_{1,1}/\partial \mbox{ }^jx=
-i\partial F_{1,2}/\partial \mbox{ }^jw$,
substituting $L^{-1}K=J$ for Equation $(2.3)$ with pair of variables
$(x,y)$ we get
$\partial F_{1,1}/\partial \mbox{ }^jx=
i\partial F_{1,1}/\partial \mbox{ }^jy$ and
$\partial F_{1,2}/\partial \mbox{ }^jx=
-i\partial F_{1,2}/\partial \mbox{ }^jy$.
Using the equality $F_{l,j}=G_{l,j}+iH_{l,j}$ we get
Equations $(2.2)$ from the latter equations.
\par Let now $F$ be differentiable at $a$ and
let $F$ be satisfying Conditions $(2.2)$.
Then $f(z)-f(a)=\sum_{l=1}^n \{
(\partial F/\partial \mbox{ }^lv)\Delta \mbox{ }^lv+$
$(\partial F/\partial \mbox{ }^lw)\Delta \mbox{ }^lw+$
$(\partial F/\partial \mbox{ }^lx)\Delta \mbox{ }^lx+$
$(\partial F/\partial \mbox{ }^ly)\Delta \mbox{ }^ly \}
+\epsilon (z-a)|z-a| $,
where $\Delta (\mbox{ }^lv,\mbox{ }^lw,\mbox{ }^lx,\mbox{ }^ly)
=\sigma ^{-1}(\mbox{ }^lz)-\sigma
(\mbox{ }^la)$ for each $l=1,...,n$.
From Conditions $(2.3)$ equivalent to $(2.2)$ we get
$f(z)-f(a)= \sum_{l=1}^n \{
(\partial F/\partial \mbox{ }^lv)I\Delta \mbox{ }^lv+$
$(\partial F/\partial \mbox{ }^lv)J\Delta \mbox{ }^lw+$
$(\partial F/\partial \mbox{ }^lv)K\Delta \mbox{ }^lx+$
$(\partial F/\partial \mbox{ }^lv)L\Delta \mbox{ }^ly \}$
$+\epsilon (z-a)|z-a| =$
$(\partial F/\partial \mbox{ }^lv)I\Delta \mbox{ }^lz$
$+\epsilon (z-a)|z-a| ,$ where $\epsilon $ is a function
continuous at $0$ and $\epsilon (0)=0$. Therefore,
$f$ is superdifferentiable by $z$ at $a$ such that
$f'(a)$ is right superlinear.
\par {\bf 2.3. Remark.} A function $f$ superdifferentiable
at each point $a\in U$ (by either $z$ or $\tilde z$ or $(z,\tilde z)$)
is called superdifferentiable
in $U$ (by either $z$ or $\tilde z$ or $(z,\tilde z)$ respectively).
The first pair of equations
in $(2.2)$ yields the usual Cauchy-Riemann conditions for
complex-valued differentiable functions. On the other hand,
we restrict $z$ by ${{t\quad 0}\choose {0\quad \bar t}},$
then $(2.1)$ also yields the usual Cauchy-Riemann condition
in complex form $D_{\bar t}f_{1,1}(t,0)=0.$
\par If a series
$f(z)=\sum_{l_1,...,l_n}\mbox{ }_{l_1}f(z)...\mbox{ }_{l_n}f(z)$
uniformly converges on $U$ together with its superdifferential, then
$$D_zf(z).h=\sum_j
\sum_{l_1,...,l_n}\mbox{ }_{l_1}f(z)...
\mbox{ }_{l_{j-1}}f(z)(D_z\mbox{ }_{l_j}f(z).h)
\mbox{ }_{l_{j+1}}f(z)...\mbox{ }_{l_n}f(z),$$
where each $\mbox{ }_lf$ is supposed to be superdifferentiable.
A similar equality holds for $D_{\tilde z}f$. This illustrates, that
in general a product of $\bf H$-valued functions need not
have a right superlinear superdifferential on the superalgebra
$\bf H^n$ even if each $\mbox{ }_lf$ has that property.
Nevertheless such functions $f$ satisfy the superdifferentiability
conditions of Definition $(2.1)$.
\par {\bf 2.4. Corollary.} {\it Let $f$ be a continuously
superdifferentiable function by $z$
with a right superlinear superdifferential on the superalgebra
$\bf H^n$ in an open subset $U$ in $\bf H^n$ and let $F$
be twice continuously differentiable by $(v,w,x,y)$ in $U$, then
certain components of $F$ are harmonic functions
by pairs of variables $(v,w)$, $(w,x)$, $(x,y)$
and $(v,y),$ namely:
$$(2.4)\quad \bigtriangleup _{\mbox{ }^lv,\mbox{ }^lw}G_{1,1}=0,\quad
\bigtriangleup _{\mbox{ }^lv,\mbox{ }^lw}H_{1,1}=0,\quad
\bigtriangleup _{\mbox{ }^lv,\mbox{ }^lw}G_{1,2}=0, \quad
\bigtriangleup _{\mbox{ }^lv,\mbox{ }^lw}H_{1,2}=0,$$
$$\bigtriangleup _{\mbox{ }^lw,\mbox{ }^lx}G_{1,1}=0,\quad
\bigtriangleup _{\mbox{ }^lw,\mbox{ }^lx}H_{1,1}=0,\quad
\bigtriangleup _{\mbox{ }^lw,\mbox{ }^lx}G_{1,2}=0, \quad
\bigtriangleup _{\mbox{ }^lw,\mbox{ }^lx}H_{1,2}=0,$$
$$\bigtriangleup _{\mbox{ }^lx,\mbox{ }^ly}G_{1,1}=0,\quad
\bigtriangleup _{\mbox{ }^lx,\mbox{ }^ly}H_{1,1}=0,\quad
\bigtriangleup _{\mbox{ }^lx,\mbox{ }^ly}G_{1,2}=0, \quad
\bigtriangleup _{\mbox{ }^lx,\mbox{ }^ly}H_{1,2}=0,$$
$$\bigtriangleup _{\mbox{ }^lv,\mbox{ }^ly}G_{1,1}=0,\quad
\bigtriangleup _{\mbox{ }^lv,\mbox{ }^ly}H_{1,1}=0,\quad
\bigtriangleup _{\mbox{ }^lv,\mbox{ }^ly}G_{1,2}=0, \quad
\bigtriangleup _{\mbox{ }^lv,\mbox{ }^ly}H_{1,2}=0$$
for each $l=1,..,n,$ where
$\bigtriangleup _{\mbox{ }^lv,\mbox{ }^lw}G_{1,1}:=
\partial ^2G_{1,1}/\partial \mbox{ }^lv^2+
\partial ^2G_{1,1}/\partial \mbox{ }^lw^2$.}
\par {\bf Proof.} From the first row of $(2.2)$
and in view of the twice continuous differentiability of $F$
it follows, that $\partial ^2G_{1,1}/\partial \mbox{ }^lv^2=
\partial ^2H_{1,1}/\partial \mbox{ }^lv\partial \mbox{ }^lw=$
$\partial ^2H_{1,1}/\partial \mbox{ }^lw\partial \mbox{ }^lv=$
$-\partial ^2G_{1,1}/\partial \mbox{ }^lw^2$.
Analgously, from the remaining rows of $(2.2)$ we deduce
the other equations in $(2.4)$. The latter
equations follow from
$\bigtriangleup _{\mbox{ }^lv,\mbox{ }^ly}=
\bigtriangleup _{\mbox{ }^lv,\mbox{ }^lw}-
\bigtriangleup _{\mbox{ }^lw,\mbox{ }^lx}+
\bigtriangleup _{\mbox{ }^lx,\mbox{ }^ly}.$
\par {\bf 2.5. Note and Definition.} Let $U$ be an open
subset in $\bf H$ and let $f: U\to \bf H$ be a function defined on
$U$ such that
$$(2.5.i)\quad f(z, {\tilde z})=
f^1(z, {\tilde z})...f^l(z, {\tilde z}),$$ where each function
$f^s(z, {\tilde z})$ is presented by a Laurent series
$$(2.5.ii)\quad f^s(z,{\tilde z})=\sum_{n=n_0}^{\infty }\sum_{m=m_0}^{\infty }
f^s_{n,m}(z-\zeta )^n(\tilde z -\tilde \zeta )^m$$
converging on $U$, where $f^s_{n,m}\in \bf H$, $z\in U$, $\zeta \in \bf H$
is a marked point, $m\in \bf Z$, $n\in \bf Z$,
if $\min (n_0,m_0)<0$, then $\zeta \notin U$. Consider the case
$f^s_{-1,m}=0$ for each $m$ and $s$. The case with terms
$f^s_{-1,m}\ne 0$ will be considered later.
\par Let $[a,b]$ be a segment in $\bf R$ and $\gamma : [a,b]\to \bf H$
be a continuous function.
Consider a partitioning $P$ of
$[a,b]$, that is, $P$ is a finite subset of $[a,b]$
consisting of an increasing sequence of points
$a=c_0<...<c_k<c_{k+1}<...<c_q=b$, then the norm of $P$
is defined as $|P|:=\max _k(x_{k+1}-x_k)$ and the $P$-variation
of $\gamma $ as $v(\gamma ;P):=\sum_{k=0}^{q-1}
|\gamma (c_{k+1})-\gamma (c_k)|,$
where $q=q(P)\in \bf N$.
The total variation (or the length) of $\gamma $ is defined as
$V(\gamma )=\sup_Pv(\gamma ;P)$. Suppose that $\gamma $ is rectifiable,
that is, $V(\gamma )<\infty $. For $f$  having decomposition
$(2.5)$ with
$f^s_{-1,m}=0$ for each $m$ and $s$
and a rectifiable path $\gamma : [a,b]\to U$
we define a (noncommutative) quaternion line integral by the formula:
$$(2.6)\quad \int_{\gamma }f(z,{\tilde z})dz:=
\lim_P I(f,\gamma ;P),\mbox{ where}$$
$$(2.7)\quad I(f,\gamma ;P):=\sum_{k=0}^{q-1}
{\hat f}(z_{k+1},{\tilde z}_{k+1}).(\Delta z_k),$$
${\hat f}(z,{\tilde z}).h:=(\partial g(z,{\tilde z})/\partial z).h$
for each $h\in \bf H$ and each $s$,
where $\Delta z_k:=z_{k+1}-z_k$, $z_k:=\gamma (c_k)$ for each
$k=0,...,q$, and where without loss of generality we suppose, that
$g$ is a function such that $(\partial g(z, {\tilde z}) / \partial z).I=
f(z,{\tilde z})$ for each $z\in U$.
In a similar way we define $\int_{\gamma }f(z,{\tilde z})d\tilde z$.
We may write shortly $\int_{\gamma }f(z)dz$ or $\int_{\gamma }
f(z)d{\tilde z}$ also for such integrals due to the bijective
correspondence between $z$ and $\tilde z$.
\par This definition is justified by the following proposition.
\par {\bf 2.6. Proposition.} {\it  Let $f$ be a function
as in \S 2.5 and suppose that there are two constants $r$ and $R$ such
that the Laurent series $(2.5)$ converges in the set
$B(a,r,R,{\bf H}):=\{ z\in {\bf H}: r\le |z-a|\le R \} $ for each
$s=1,...,l$, let also $\gamma $ be a rectifiable path contained in
$U\cap B(a,r',R',{\bf H})$, where $r<r'<R'<R$.
Then the quaternion line integral exists.}
\par {\bf Proof.} Since each $f^s$ converges in $B(a,r,R,{\bf H})$, then
$${\overline {\lim }}_{n>0, m>0}|f^s_{n,m}|^{1/(n+m)}R\le 1,\mbox{ hence}$$
$$\| f \|_{\omega }:=\prod_{s=1}^l (\sup_{n+m<0}|f^s_{n,m}| r^{n+m},
\sup_{n+m\ge 0} |f^s_{n,m}| R^{n+m} )<\infty $$
 and inevitably
$$\| f \| _{1,\omega ,B(a,r',R',{\bf H})}
:=\prod_{s=1}^l [(\sum_{n+m<0}|f^s_{n,m}|{r'}^{n+m})
+(\sum_{n+m>0}|f^s_{n,m}| {R'}^{n+m})] <\infty .$$
For each locally $(z,{\tilde z})$-analytic
function $f$ in $U$ and each $z_0$ in $U$ there exists a ball of radius
$r>0$ with center $z_0$ such that $f$ has a decomposition
analogous to $(2.5i,ii)$ in this ball with all $n_j$ and $m_j$ nonnegative,
$j=1,..,l$. Consider two quaternion $(z,\tilde z)$-locally analytic
functions $f$ and $q$ on $U$ such that $f$ and $q$ noncommute.
Let $f^0:=f$, $q^0:=q$, $q^{-n}:=q^{(n)}$, $\partial (q^n)/\partial z=:
q^{n-1}$ and $q^{-k-1}=0$ for some $k\in \bf N$, then
\par $(i)$ $(fq)^1=f^1q-f^2q^{-1}+f^3q^{-2}+...+(-1)^kf^{k+1}q^{-k}$.
In particular, if $f=az^n$, $q=bz^k$, with $n>0$, $k>0$,
$b\in {\bf H}\setminus {\bf R}I$, then $f^p=[(n+1)...(n+p)]^{-1}az^{n+p}$
for each $p\in \bf N$, $q^s=(k-1)...(k-s+1)bz^{k-s}$.
Also
\par $(ii)$ $(fq)^1=fq^1-f^{-1}q^2+f^{-2}q^3+...+(-1)^pf^{-p}q^{p+1}$.
Apply $(i)$ for $n\ge m$ and $(ii)$ for $n<m$
to solve the equation $(\partial g(z, {\tilde z})/\partial z).I=
f(z,{\tilde z})$ for each $z\in U$. If $f$ and $q$
have series converging in $Int (B(z_0,r,{\bf H}))$,
then these formulas show that there exists a $(z,\tilde z)$-analytic
function $(fq)^1$ with series converging in $Int (B(0,r,{\bf H}))$,
since $\lim_{n\to \infty }(nr^n)^{1/n}=r$, where $0<r<\infty $. 
Consider the equation $BA=AC$, where $A$, $B$ and $C$ are quaternions.
Therefore, for each quaternion locally $(z,\tilde z)$-analytic function
$f$ there exists the operator $\hat f$. Considering a function
$G$ of real variables corresponding to $g$ we get that all solutions
$g$ differ on quaternion constants, hence $\hat f$ is unique for $f$.
If $A\ne 0$, then $C=A^{-1}BA$, hence $|C|=|B|$. If $B\ne 0$, then
$C=BD$, where $D=B^{-1}C$ and $|D|=1$.
Therefore,
$$(2.8)\quad f^s_{n,m}(z_{j+1}-a)^k(\Delta z_j)(z_{j+1}-a)^{n-k}
(\tilde z_{j+1}-\tilde a)^m=$$
$$f^s_{n,m}(z_{j+1}-a)^n(\tilde z_{j+1}-\tilde a)^mC(n-k,m;z_{j+1},z_j,a)
(\Delta z_j),$$ where $C(p,m;z_{j+1},z_j,a)\in \bf H$ and
$|C(p,m;z_{j+1},z_j,a)|=1$ for each $z_{j+1}\ne z_j$, $z_{j+1}\ne a$,
for each $p$ and $m$.
From Equation $(2.8)$ it follows, that 
$|I(f,\gamma ;P)| \le \| f \| _{1,\omega ,B(a,r',R',{\bf H})} v(\gamma ;P),$
for each $P$, and inevitably
$$(2.9)\quad |I(f,\gamma ;P)-I(f,\gamma ;Q)|\le w(\hat f;P)V(\gamma )$$
for each $Q\supset P$, where
$$(2.10)\quad w(\hat f;P):=\max_{(z, \zeta \in \gamma ([c_j,c_{j+1}]))}
\{ \| \hat f(z,{\tilde z})-\hat f(\zeta ,{\tilde \zeta }) \|: \quad
z_j=\gamma (c_j), c_j\in P \} ,$$
$ \| {\hat f}(z,{\tilde z}) -{\hat f}(\zeta ,{\tilde \zeta }) \|
:= \sup_{h\ne 0}|{\hat f}(z,{\tilde z}).h
-{\hat f}(\zeta ,{\tilde \zeta }).h |/|h|$. Since $\lim_{n\to \infty }
(n)^{1/n}=1,$ then $\lim_P \omega ({\hat f},P)=0.$
From $\lim_Pw(\hat f;P)=0$ the existence of $\lim_PI(f,\gamma ;P)$
now follows.
\par {\bf 2.7. Theorem.} {\it Let $\gamma $ be a rectifiable
path in $U$, then the quaternion line integral
has a continuous extension on the space $C^0_b(U,{\bf H})$
of bounded continuous functions $f: U \to \bf H$.
This integral is an $\bf R$-linear and left-$\bf H$-linear
functional on $C^0_b(U,{\bf H})$.}
\par {\bf Proof.} Since $\gamma $ is continuous and $[a,b]$ is compact,
then there exists a compact canonical closed subset $V$ in $\bf H$, that is,
$cl (Int(V))=V$, such that $\gamma ([a,b])\subset V\subset U$.
Let $f\in C^0_b(U,{\bf H})$,
then in view of the Stone-Weierstrass theorem
for a function $F(v,w,x,y)=f\circ \sigma (v,w,x,y)$ and each
$\delta >0$ there exists a polynomial $T$ such that $\| F- T\|_0<\delta $,
where $\| f\|_0:=\sup_{z\in U}|f(z)|$. This polynomial takes values in
$\bf H$, hence it has the form:
$T_{1,1}={\bar T}_{2,2}=\alpha ^{1,1}_{i_1,i_2,i_3,i_4}v^{i_1}w^{i_2}
x^{i_3}y^{i_4}$ and
$T_{1,2}=-{\bar T}_{2,1}=\alpha ^{1,2}_{i_1,i_2,i_3,i_4}v^{i_1}w^{i_2}
x^{i_3}y^{i_4}$,
where summation is accomplished by repeated upper and lower
indices. There are relations $J^2=K^2=L^2=-I$, $JK=-KJ=L$,
$KL=-LK=J$, $LJ=-JL=K$, consequently, $zJ=vJ-wI-xL+yK,$
$zK=vK+wL-xI-yJ,$ $zL=vL-wK+xJ-yI,$ where $z=vI+wJ+xK+yL$
is in $\bf H$, $v$,$w$, $x$ and $y$ are in $\bf R$.
Therefore, $\tilde z=vI-wJ-xK-yL$,
$J\tilde z=vJ+wI-xL+yK$, $K\tilde z=vK+wL+xI-yJ,$
$L\tilde z=vL-wK+xJ+yI$, hence $(z+\tilde z)/2=vI,$
$(J\tilde z-zJ)/2=wI,$ $(K\tilde z-zK)/2=xI$,
$(L\tilde z-zL)/2=yI$. From this it follows, that
$T$ can be expressed in $z$ and $\tilde z$
such that
$$T=[Re (\alpha ^{1,1}_{i_1,i_2,i_3,i_4})I+
Im (\alpha ^{1,1}_{i_1,i_2,i_3,i_4})J +
Re (\alpha ^{1,2}_{i_1,i_2,i_3,i_4})K+
Im (\alpha ^{1,1}_{i_1,i_2,i_3,i_4})L]$$
$$\times [(z+\tilde z)/2]^{i_1}[(J\tilde z-zJ)/2]^{i_2}
[(K\tilde z-zK)/2]^{i_3}[(L\tilde z-zL)/2]^{i_4}.$$
This polynomial can be rewritten in a form similar to $f$ in \S 2.5
(see Formulas $(2.5.i,ii)$), since $z$ and $\tilde z$ commute.
Two variables $(z+\tilde z)/2=vI$ and $(z-\tilde z)/2=wJ+xK+yL$
commute for each $z\in \bf H$.
Therefore, the $\bf R$-linear space of functions
on $U$ having decomposition $(2.5.i,ii)$ is dense in $C^0_b(U,{\bf H})$.
\par Consider a function $g(z,\tilde z)$ on $U$, suppose
that $q(z,\tilde \zeta)$ is another function on $U^2$ such that
$q(z,\tilde \zeta )|_{z=\zeta }=g(z,\tilde z)$.
Let $q(z,\tilde \zeta )$ be superdifferentiable by $z$
for a fixed $\zeta \in U$, then $\partial q(z,\tilde \zeta )/\partial z$
for $z=\zeta $ is denoted by $\partial g(z,\tilde z)/\partial z$.
Consider a space of all such that $g$ on $U$ for which
$(\partial g(z,\tilde z)/\partial z).S$ is a bounded continuous function
on $U$ for each $S\in \{ I,J,K,L \} $,
it is denoted by $C^1_b(U,{\bf H})=C^{1,0}_b(U,{\bf H})$
and it is supplied with the norm $\| g \| _{C^1_b}:=\| g\|_{C^0_b}+
\sum_{S\in \{ I,J,K,L \} }
\| (\partial g(z,\tilde z)/\partial z).S \|_{C^0_b}$, where
$\| g \|_{C^0_b}:=\sup_{z\in U}|g(z)|$.
In view of Proposition 2.2 for each $g\in C^1_b(U,{\bf H})$
the corresponding function $Q(z,\zeta )$ satisfies Condition
$(2.1)$ by $z$. This entails, that $\partial Q/\partial v$,
$\partial Q/\partial w$, $\partial Q/\partial x$
and $\partial Q/\partial y$ are in $C^0_b(U^2,{\bf H})$ (see \S 2.1).
Consequently, imposing the condition $z=\zeta :$ $\quad $
$(\partial g/\partial z).J$, $(\partial g/\partial z).K$
and $(\partial g/\partial z).L$ are also continuous bounded functions,
hence $(\partial g/\partial z).h\in C^0_b(U\times B(0,0,1,{\bf H}),{\bf H}),$
where $h\in B(0,0,1,{\bf H})$. Therefore, there exists
a positive constant $C$ such that
$$(2.11)\quad \sup_{h\ne 0} | (\partial g/\partial z).h|/|h| \le C
\sum_{S\in \{ I,J,K,L \} }
\| (\partial g/\partial z).S \|_{C^0_b} ,$$
since $h=v_hI+w_hJ+x_hK+y_hL$ for each $h\in \bf H$
and $(\partial g/\partial z)$ is $\bf R$-linear and
$(\partial g/\partial z).(h_1+h_2)=(\partial g/\partial z).h_1+
(\partial g/\partial z).h_2$ for each $h_1$ and $h_2\in \bf H$, where
$v_h$, $w_h$, $x_h$ and $y_h$ are real numbers,
$G(v,w,x,y):=g\circ \sigma (v,w,x,y)$ is Frech\'et differentiable
on an open subset $U_{\sigma }\subset \bf R^4$ such that
$\sigma (U_{\sigma })=U$.
\par In \S 2.6 it was shown that the equation
$(\partial g(z,\tilde z)/\partial z).I=f(z,\tilde z)$ has a solution
in a class of quaternion locally $(z,\tilde z)$-analytic functions on $U$.
The subset $C^{\omega }_{(z,\tilde z)}(U,{\bf H})$ is dense in
the uniform space $C^0_b(U,{\bf H})$.
\par If $g=g^1...g^l$ is a product of functions
$g^s\in C^{1,0}_b(U,{\bf H})$, then $(\partial g/\partial z).h=
\sum_{s=1}^lg^1(z,\tilde z)...g^{s-1}(z,\tilde z)
[(\partial g^s/\partial z).h]g^{s+1}(z,\tilde z)...g^l(z,\tilde z)$
for each $h\in \bf H$. 
Consider the space ${\hat C}^0_b(U,{\bf H}):= \{
((\partial g/\partial z).I,(\partial g/\partial z).J,
(\partial g/\partial z).K, (\partial g/\partial z).L):$
$g\in C^{1,0}_b(U,{\bf H}) \} $.
It has an embedding $\xi $ into $C^0_b(U,{\bf H})$ and
$\| g \| _{C^{1,0}_b}\ge \sum_{S\in \{ I,J,K,L \} }
\| (\partial g/\partial z).S \| _{C^0_b} $.
In view of Inequality $(2.11)$ the completion of
${\hat C}^0_b(U,{\bf H})$ relative to $\| * \|_{C^0_b(U,{\bf H})}$
coinsides with $C^0_b(U,{\bf H})$.
\par Let $\{ f^j: \quad j\in {\bf N} \} $ be a sequence of functions
having decomposition $(2.5)$ and converging to $f$ in
$C^0_b(U,{\bf H})$ relative to the metric $\rho (f,q):=\sup_{z\in U}
|f(z,{\tilde z})-q(z,{\tilde z})|$ such that $f^j=\xi 
((\partial g^j/\partial z).I,(\partial g^j/\partial z).J,
(\partial g^j/\partial z).K, (\partial g^j/\partial z).L)$
for some $g^j\in C^{1,0}_b(U,{\bf H})$.
Relative to this metric $C^0_b(U,{\bf H})$ is complete.
We have the equality
$$\partial (\int_0^sF(v_0+\phi h_v,w_0+\phi h_w,x_0+\phi h_x,
y_0+\phi h_y)d\phi )/\partial s =F(v,w,x,y)$$
for each continuous
function $F$ on $U_{\sigma }$, where $v=v_0+sh_v$, $w=w_0+sh_w$,
$x=x_0+sh_x$ and $y=y_0+sh_y$, $(v_0,w_0,x_0,y_0)+\phi (h_v,h_w,h_x,h_y)
\in U_{\sigma }$ for each $\phi \in \bf R$ with $0\le \phi \le s$,
$h_v$, $h_w$, $h_x$ and $h_y\in \bf R^4$.
Let $z_0$ be a marked point in $V$.
There exists $R>0$ such that $\gamma $ is contained in the
interior of the parallelepiped
$V:=\{ z\in {\bf H}: z=vI+wJ+xK+yL, |v-v_0|\le R,
|w-w_0|\le R, |x-x_0|\le R, |y-y_0|\le R \} $.
\par If $V$ is not contained in
$U$ consider a continuous extension of a continuous function $F$
from $V\cap U_0$ on $V$, where $U_0$ is a closed subset in $U$
such that $Int (U_0)\supset \gamma $ (about the theorem
of a continuous extension see \cite{eng}).
Therefore, suppose that $F$ is given on $V$. Then the function
$F_1(v,w,x,y):=\int_{v_0}^v\int_{w_0}^w\int_{x_0}^x\int_{y_0}^y
F(v_1,w_1,x_1,y_1)dv_1dw_1dx_1dy_1$ is in $C^1(V,{\bf H})$
(with one sided derivatives on $\partial V$ from inside $V$).
Consider a foliation of $V$ by three dimensional
$C^0$-manifolds $\Upsilon _z$ such that $\Upsilon _z\cap \Upsilon _{z_1}=
\emptyset $ for each $z\ne z_1$, where $z, z_1\in \gamma $,
$\bigcup_{z\in \gamma }\Upsilon _z=V_1$, $V_1$ is a canonical closed
subset in $\bf H$ such that $\gamma \subset V_1\subset V$.
Choose this foliation such that to have decomposition
of a Lebesgue measure $dV$ into the product of measures $d\nu (z) $
along $\gamma $ and $d\Upsilon _z$ for each $z\in \gamma $.
In view of the Fubini theorem there exists $\int_Vf(v_1,...,y_1)dV=
\int_{\gamma }(\int_{\Upsilon _z}f(z,{\tilde z})d\Upsilon _z)d\nu (z)$.
If $\gamma $ is a straight line segment then
$\int_{\gamma } f(z,{\tilde z})dz$
is in $L^1(\Upsilon ,{\bf H})$.
Let $U_{\bf R}$ be a real region in $\bf R^4$ corresponding to
$U$ in $\bf H$.
\par Consider the Sobolev space $W^s_2(U_{\bf R},{\bf R^4})$
of functions $h: U_{\bf R}\to \bf R^4$ for which $D^{\alpha }h\in
L^2(U_{\bf R},{\bf R^4})$ for each $|\alpha |\le s$, where
$0\le s\in \bf Z$. In view of Theorem 18.1.24
\cite{hoermpd} (see also the notation there)
if $A\in \Psi ^m$ is a properly supported
pseudodifferential elliptic operator of order $m$ in the sence that
the principal symbol $a\in S^m(T^*(X))/S^{m-1}(T^*(X))$ has an inverse
in $S^{-m}(T^*(X))/S^{-m-1}(T^*(X))$, then one can find
$B\in \Psi ^{-m}$ properly supported such that $BA-I\in \Psi ^{-\infty }$,
$AB-I\in \Psi ^{-\infty }$. One calls $B$ a parametrix for $A$.
In view of Proposition 18.1.21 \cite{hoermpd} each $A\in \Psi ^m$
can be written as a sum $A=A_1+A_0$, where $A_1\in \Psi ^m$ is properly
supported and the kernel of $A_0$ is in $C^{\infty }$.
In particular we can take a pseudodifferential operator
with the principal symbol $a(x,\xi )=(b+|\xi |^2)^{s/2}$,
where $b>0$ is a constant and $s\in \bf Z$, which corresponds
to $b+\Delta $ for $s=1$ up to minor terms,
where $\Delta =\nabla ^2$ is the Laplacian
(see also Theorem 3.2.13 \cite{grubb} about its parametrix family).
For estimates of a solution there may be also applied Theorem 3.3.2
and Corollary 3.3.3 \cite{grubb} concerning parabolic
pseudodifferential equations for our particular case
corresponding to $(\partial g/\partial z).I=f$ rewritten in real variables.
\par Due to the Sobolev theorem (see \cite{stein,trieb})
there exists an embedding of the Sobolev space $W^3_2(V,{\bf H})$
into $C^0(V,{\bf H})$ such that $\| g \|_{C^0}\le C\| g \|_{W^3_2}$
for each $g\in W^3_2$, where $C$ is a positive constant independent of $g$.
If $h\in W^{k+1}_2(V,{\bf H})$, then $\partial h/\partial b_j
\in W^k_2(V,{\bf H})$ for each $k\in \bf N$ and in particular for $k=3$
and each $j=1,...,4$ (see \cite{stein}).
On the other hand $\| h \|_{L^2(V,{\bf H})}\le
\| h \|_{C^0(V,{\bf H})}(2R)^2$ for each $h\in L^2(V,{\bf H})$.
Therefore, $\| A^{-k}h \| _{W^k_2(V,{\bf H})}
\le C \| h \| _{C^0(V,{\bf H})}(2R)^{k+2}$ for each $k\in \bf N$,
where $C=const >0$, $A$ is an elliptic pseudodifferential
operator such that $A^2$ corresponds to $(1+\Delta )$.
From Equations $(2.6,2.7)$ and Inequality $(2.11)$ it follows, that
$$(2.12)\quad |I(f-q,\gamma ;P)|\le \rho (f,q)V(\gamma )C_1 \exp (C_2 R^6)$$
for each partitioning $P$, where $C_1$ and $C_2$ are positive
constants independent
of $R$, $f$ and $q$. In view of Formulas $(2.9,2.10)$ $\{ \int_{\gamma }
f^j(z,{\tilde z})dz: j\in {\bf N} \} $ is a Cauchy sequence in $\bf H$
and the latter is complete as the metric space. Therefore,
there exists $\lim_j\lim_PI(f^j,\gamma ;P)=\lim_j\int_{\gamma }f^j(z,
{\tilde z})dz$, which we denote by $\int_{\gamma }f(z,{\tilde z})dz$.
As in \S 2.6 we get that all solutions $g$ differ on quaternion
constants on each connected component of $U$, consequently,
the functional $\int_{\gamma }$ is uniquely defined on $C^0_b(U,{\bf H})$.
The functional $\int_{\gamma }: C^0_b(U,{\bf H})\to \bf H$
is continuous due to Formula $(2.12)$ and evidently
it is $\bf R$-linear, since $\lambda z=z\lambda $
for each $\lambda \in \bf R$ and each $z\in \bf H$,
that is, $\int_{\gamma }(\lambda _1f_1(z,{\tilde z})+
\lambda _2f_2(z,{\tilde z}))dz=$
$\int_{\gamma }(f_1(z,{\tilde z})\lambda _1+
f_2(z,{\tilde z})\lambda _2)dz=$
$\lambda _1\int_{\gamma }f_1(z,{\tilde z})dz+
\lambda _2\int_{\gamma }f_2(z,{\tilde z})dz$
for each $\lambda _1$ and $\lambda _2\in \bf R$,
$f_1$ and $f_2\in C^0_b(U,{\bf H})$.
Moreover, it is left-$\bf H$-linear,
that is, $\int_{\gamma }(\lambda _1f_1(z,{\tilde z})+
\lambda _2f_2(z,{\tilde z}))dz=
\lambda _1\int_{\gamma }f_1(z,{\tilde z})dz+
\lambda _2\int_{\gamma }f_2(z,{\tilde z})dz$
for each $\lambda _1$ and $\lambda _2\in \bf H$,
$f_1$ and $f_2\in C^0_b(U,{\bf H})$, since $I(f,\gamma ;P)$ is
left-$\bf H$-linear.
\par {\bf 2.8. Remark.} Let $\eta $ be a differential form
on open subset $U$ of the Euclidean space $\bf R^{4m}$ with values in
$\bf H$, then it can be written as
$$(2.13)\quad \eta =\sum_{\Upsilon }\eta _{\Upsilon }db^{\wedge \Upsilon },$$
where $b=(\mbox{ }^1b,...,\mbox{ }^mb)\in \bf R^{4m}$,
$\mbox{ }^jb=(\mbox{ }^jb_1,...,\mbox{ }^jb_4)$,
$\mbox{ }^jb_i\in \bf R$, $\eta _{\Upsilon }=\eta _{\Upsilon } (b):
{\bf R^{4m}}\to \bf H$ are $s$ times continuously differentiable
$\bf H$-valued functions with $s\in \bf N$,
$\Upsilon =(\Upsilon (1),...,\Upsilon (m))$, $\Upsilon (j)=
(\Upsilon (j,1),...,\Upsilon (j,4))\in
\bf N^4$ for each $j$, $db^{\wedge \Upsilon }=
d\mbox{ }^1b^{\wedge \Upsilon (1)}\wedge ... \wedge d\mbox{ }^m
b^{\wedge \Upsilon (m)}$, $d\mbox{ }^jb^{\wedge \Upsilon (j)}=
d\mbox{ }^jb_1^{\Upsilon (j,1)}\wedge ... \wedge d\mbox{ }^j
b_4^{\Upsilon (j,4)}$, where $d\mbox{ }^jb_i^0=1$, 
$d\mbox{ }^jb_i^1=d\mbox{ }^jb_i$, $d\mbox{ }^jb_i^k=0$
for each $k>1$. If $s\ge 1$, then there is defined
an (external) differential
$$d\eta =\sum_{\Upsilon ,(j,i)}(\partial \eta _{\Upsilon }/\partial
\mbox{ }^jb_i)(-1)^{\alpha (j,i)}db^{\wedge (\Upsilon +e(j,i))} ,$$
where $e(j,i)=(0,...,0,1,0,...,0)$ with $1$ on the
$4(j-1)+i$-th place, $\alpha (j,i)=(\sum_{l=1}^{j-1}\sum_{k=1}^4
\Upsilon (l,k))+\sum_{k=1}^{i-1}\Upsilon (j,k)$.
These differential forms have matrix structure themselves.
Consider basic matrices $S=(\mbox{ }^1S,...,\mbox{ }^mS)$
and their ordered product $S^{\to \Upsilon }:=
\mbox{ }^1S^{\to \Upsilon (1)}...\mbox{ }^mS^{\to \Upsilon (m)}$,
where $\mbox{ }^jS=(\mbox{ }^jS_1,...,\mbox{ }^jS_4)=
(I,J,K,L)$, $\mbox{ }^jS^{\to \Upsilon (j)}=J^{\Upsilon (j,2)}
K^{\Upsilon (j,3)}L^{\Upsilon (j,4)}$, $S^0=I$.
Then Equation $(2.13)$ can be rewritten in the form:
$$(2.14)\quad \eta =\sum_{\Upsilon }\xi _{\Upsilon }
d(Sb)^{\wedge \Upsilon },$$
where $Sb=(\mbox{ }^1S_1\mbox{ }^1b_1,...,\mbox{ }^mS_m
\mbox{ }^mb_m)\in \bf H^{4m}$, $d\mbox{ }^jS_k\mbox{ }^jb_k=
\mbox{ }^jS_kd\mbox{ }^jb_k$, $\xi _{\Upsilon }:=
\eta _{\Upsilon }(S^{\to \Upsilon })^{-1}$.
Relative to the external product $Id\mbox{ }^jb_1$
anticommutes with others basic differential $1$-forms
$\mbox{ }^jS_kd\mbox{ }^jb_k$; for $k=2,3,4$
these forms commute with each other relative to the external product.
This means that the algebra of quaternion differential forms is graded
relative to the external product.
\par From \S 2.7 it follows, that $(dz+d\tilde z)/2=Idv,$
$Jdw=-(d\tilde z+J(dz)J)/2,$ $Kdx=-(d\tilde z+K(dz)K)/2$,
$Ldy=-(d\tilde z+L(dz)L)/2$. Therefore, the right side of
Equation $(2.14)$ can be rewritten with $d\mbox{ }^jz$,
$d\mbox{ }^j\tilde z$,
$Jd\mbox{ }^jzJ$, $Kd\mbox{ }^jzK$ and $Ld\mbox{ }^jzL$
on  the right side. From the
latter 5 differential $1$-forms 4 linearly independent ones can be
chosen, since summing these forms we have:
$dz=-2 d\tilde z-JdzJ-KdzK-LdzL$, hence $LdzL=-(dz+2 d\tilde z+
JdzJ+KdzK)$. These $1$-forms do neither commute nor anticommute,
since they are not pure elements of the graded algebra.
For example,
$(dz\wedge dz)_{1,1}={\overline {(dz\wedge dz)}}_{2,2}=
-du\wedge d\bar u$, $(dz\wedge dz)_{1,2}=-
{\overline {(dz\wedge dz)}}_{2,1}=dt\wedge du-d\bar t\wedge du$;
$(dz\wedge d\tilde z)_{1,1}={\overline
{(dz\wedge d\tilde z)}}_{2,2}=dt\wedge d\bar t
+du\wedge d\bar u$, $(dz\wedge d\tilde z)_{1,2}=-
{\overline {(dz\wedge d\tilde z)}}_{2,1}=-2dt\wedge du$;
$(d\tilde z\wedge dz)_{1,1}={\overline
{(d\tilde z\wedge dz)}}_{2,2}=d\bar t\wedge dt
+du\wedge d\bar u$, $(d\tilde z\wedge dz)_{1,2}=-
{\overline {(d\tilde z\wedge dz)}}_{2,1}=2d\bar t\wedge du$;
$(dz\wedge JdzJ)_{1,1}={\overline
{(dz\wedge JdzJ)}}_{2,2}=-du\wedge d\bar u$,
$(dz\wedge JdzJ)_{1,2}=-
{\overline {(dz\wedge JdzJ)}}_{2,1}=dt\wedge du+d\bar t\wedge du$;
$(JdzJ\wedge dz)_{1,1}={\overline
{(JdzJ\wedge dz)}}_{2,2}=-du\wedge d\bar u$,
$(JdzJ\wedge dz)_{1,2}=-
{\overline {(JdzJ\wedge dz)}}_{2,1}=-dt\wedge du-
d\bar t\wedge du$;
$(dz\wedge KdzK)_{1,1}={\overline
{(dz\wedge KdzK)}}_{2,2}=-dt\wedge d\bar t$,
$(dz\wedge KdzK)_{1,2}=-{\overline {(dz\wedge KdzK)}}_{2,1}
=-dt\wedge d\bar u+dt\wedge du$;
$(KdzK\wedge dz)_{1,1}={\overline {(KdzK\wedge dz)}}_{2,2}=
dt\wedge d\bar t$,
$(KdzK\wedge dz)_{1,2}=-{\overline {(KdzK\wedge dz)}}_{2,1}=
-d\bar t\wedge du+d\bar t\wedge d\bar u$.
On the other hand Equation $(2.13)$ can be rewritten
using the identities: $(dz+d\tilde z)/2=Idv$,
$Idw=(Jd\tilde z-dzJ)/2$, $Idx=(Kd\tilde z -dz K)/2$,
$Idy=(Ld\tilde z-dzL)/2$, where $dz={{dt \quad du}\choose
{-d{\bar u}\quad d{\bar t}}}$ in $2\times 2$ complex
matrix notation.
\par Consider a $C^1$-function $f$ on $U$ with values in $\bf H$,
then
$$D_{A\mbox{ }^j\tilde z}f=(\partial f/\partial (A\mbox{ }^j\tilde z)).
A\mbox{ }^j\tilde z=(\partial f/\partial \mbox{ }^j\tilde z).(A^{-1}A
\mbox{ }^j\tilde z)=D_{\mbox{ }^j\tilde z}f$$
for each $A\in \bf H$, also
$$D_{\mbox{ }^jzA}f=(\partial f/\partial (\mbox{ }^jzA)).
\mbox{ }^jzA=(\partial f/\partial \mbox{ }^jz).(
\mbox{ }^jz)AA^{-1}=D_{\mbox{ }^jz}f$$
for each $A\in \bf H$, where $D_hf:=(\partial f/\partial h).h$.
We apply this also in particular to $Iv$, $wJ$, $xK$ and $yL$.
\par There is the standard embedding of
the algebra of complex $n\times n$ matrices $A$ into the algebra
of real $2n\times 2n$ matrices $B$ such that in its block form
$B_{1,1}=Re(A)$, $B_{2,2}=Re (A)$, $B_{1,2}=Im (A),$
$B_{2,1}=-Im (A)$, where $B_{i,j}$ are $n\times n$ blocks.
Therefore, quaternion differential forms can be
embedded into the algebra of differential forms
over the algebra of real $4\times 4$ matrices.
This shows, that the exterior differentiation operator
$\mbox{ }_{\bf H}d$ for $\bf H$-valued differential forms
over $\bf H$ and that of
for their real matrix realization $\mbox{ }_{\bf R}d$
coincide and their common operator is denoted by $d$.
Consider the equality
$$(\partial \eta _{\Upsilon }/\partial \mbox{ }^jb^l)
\mbox{ }^jb^l \wedge db^{\Upsilon }=
[(\partial \eta _{\Upsilon }/\partial \mbox{ }^jz).
(\partial \mbox{ }^jz/\partial \mbox{ }^jb^l)]
\mbox{ }^jb^l \wedge db^{\Upsilon }$$ $$+
[(\partial \eta _{\Upsilon }/\partial \mbox{ }^j{\tilde z}).
(\partial \mbox{ }^j{\tilde z}/\partial \mbox{ }^jb_l)]
\mbox{ }^jb^l \wedge db^{\Upsilon }.$$
Applying it to $l=1,...,4$ and summing left anf right
parts of these equalities we get
$d\eta (z,\tilde z)=((\partial \eta /\partial z).
d\mbox{ }^jz)\wedge db^{\Upsilon }+ ((\partial \eta /\partial {\tilde z}).
d\mbox{ }^j{\tilde z})\wedge db^{\Upsilon }$,
hence the external differentiation can be presented in the form
$$(2.15)\quad d=\partial _z+\partial _{\tilde z},$$ where $\partial _z$
and $\partial _{\tilde z}$ are external differentiations
by variables $z$ and $\tilde z$ respectively.
\par Certainly for an external product $\eta _1\wedge \eta _2$
there is not (in general) an $\lambda \in \bf H$
such that $\lambda \eta _2\wedge \eta _1=
\eta _1\wedge \eta _2$, if $\eta _1$ and $\eta _2$ are not
pure elements (even or odd) of the graded algebra.
\par {\bf 2.9. Definition.} A Hausdorff topological space
$X$ is said to be $n$-connected for $n\ge 0$
if each continuous map $f: S^k\to X$ from the $k$-dimensional
real unit sphere into $X$ has a continuous extension
over $\bf R^{k+1}$ for each $k\le n$. A $1$-connected space
is also said to be simply connected.
\par {\bf 2.10. Remark.} In accordance with Theorem
1.6.7 \cite{span} a space $X$ is $n$-connected
if and only if it is path connected and $\pi _k(X,x)$ is trivial
for every base point $x\in X$ and each $k$ such that $1\le k\le n$.
\par Denote by $Int (U)$ an interior of a subset $U$
in a topological space $X$, by $cl (U)=\bar U$ a closure
of $U$ in $X$. For a subset $U$ in $\bf H$, let
$\pi _{1,t}(U):=\{ u: z\in U, \mbox{ where }z_{1,1}={\bar z}_{2,2}
=t, z_{1,2}=-{\bar z}_{2,1}=u \} $ for a given $t\in \bf C$;
$\pi _{2,u}(U):=\{ t: z\in U, \mbox{ where }z_{1,1}={\bar z}_{2,2}
=t, z_{1,2}=-{\bar z}_{2,1}=u \} $ for a given $u\in \bf C$,
that is, geometrically
$\pi _{1,t}(U)$ and $\pi _{2,u}(U)$ are projections on
complex planes $\bf C_2$ and $\bf C_1$ of intersections
of $U$ with planes ${\tilde \pi }_{1,t}\ni {{t \quad 0}\choose
{0\quad {\bar t}}}$ and
${\tilde \pi }_{2,u} \ni {{0 \quad u}\choose {-{\bar u}\quad 0}}$
parallel to $\bf C_2$ and $\bf C_1$ respectively.
\par {\bf 2.11. Theorem.} {\it Let $U$ be a domain in $\bf H$
such that $\emptyset \ne Int (U)\subset U\subset cl (Int (U))$
and $U$ is $3$-connected; $\pi _{1,t}(U)$ and $\pi _{2,u}(U)$
are simply connected in $\bf C$ for each
$t$ and $u\in \bf C$ for which there exists $z\in U$.
Let $f$ be a bounded continuous function
from $U$ into $\bf H$ which satisfies Condition $(2.1)$
on an open domain $W$ such that $W\supset U$. Then for each
rectifiable closed path $\gamma $ in $U$ a quaternion
line integral $\int_{\gamma }f(z)dz=0$ is equal to zero.}
\par {\bf Proof.} For a path $\gamma $ there exists
a compact canonical closed subset in $\bf H$: $\quad W\subset Int (U)$
such that $\gamma ([0,1])\subset W$, since $\gamma $
is rectifiable and $\bf H$ is locally compact.
In view of Theorem 2.7 for each sequence of functions
$f_n\in C^{1,1}(U,{\bf H})$ converging to $f$ in $C^0_b(U,{\bf H})$
such that $f_n(z)=(\partial g_n(z,\tilde z)/\partial \tilde z).I$
with $g_n\in C^{2,1}(U,{\bf H})$
and each sequence of paths $\gamma _n: [0,1] \to U$
$C^3$-continuously differentiable and converging to $\gamma $
relative to the total variation $V(\gamma -\gamma _n)$
there exists $\lim_n\int_{\gamma _n}f_n(z)dz=\int_{\gamma }f(z)dz$.
Therefore, it is sufficient to consider the case of $f\in C^{1,1}(U,
{\bf H})$ such that $f(z)=(\partial g(z,\tilde z)/\partial {\tilde z}).I$
on $U$, and continuously differentiable $\gamma $.
Denote the integral $\int_{\gamma }f(z)dz$ by $Q$,
then $Q=0$ if and only if $Q\tilde Q=0$. On  the other hand,
$\tilde Q=\lim_P{\tilde I}(f,\gamma ;P)$
and $Q=\int_{\gamma }(\partial g(z,\tilde z)/\partial z).dz$,
hence $\tilde Q=\int_{\gamma }d\tilde z
(\partial _L {\tilde g} (z,\tilde z)/\partial {\tilde z})=
\int_{\tilde \gamma }dz (\partial _L {\tilde g} (z,\tilde z)/\partial z)$,
where $(\partial _L q(z,\tilde z)/\partial z)$ is the left
derivative, ${\tilde q}(z,\tilde z):=\tilde a$, where $a=q(z,\tilde z)$.
We can write this integral in the form $Q=\int_0^1
(\partial g(z,\tilde z)/\partial z).\gamma '(t)dt$.
Denoting components of $(\partial g(z,\tilde z)/\partial z)$
as the complex $2\times 2$ matrix with entries ${\hat f}_{i,j}$
with $i$ and $j\in \{ 1,2 \} $, we get
$$Q_{1,1}={\bar Q}_{2,2}=\int_0^1{\hat f}_{1,1}(t(\theta ),u(\theta ))
{\gamma '}_{1,1}(\theta )d\theta -$$ 
$$\int_0^1{\hat f}_{1,2}(t(\theta ),u(\theta ))
{{\bar \gamma }'}_{1,2}(\theta )d\theta ;$$
$$Q_{1,2}=-{\bar Q}_{2,1}=\int_0^1{\hat f}_{1,1}(t(\theta ),u(\theta ))
{\gamma '}_{1,2}(\theta )d\theta +$$ 
$$\int_0^1{\hat f}_{1,2}(t(\theta ),u(\theta ))
{{\bar \gamma }'}_{1,1}(\theta )d\theta ,$$ where $t(\theta )=
\gamma _{1,1}(\theta )$, $u(\theta )=\gamma _{1,2}(\theta )$.
Evidently $\gamma (1)=\gamma (0)$ if and only if two equalities are
satisfied $\gamma _{1,1}(1)=\gamma _{1,1}(0)$ and
$\gamma _{1,2}(1)=\gamma _{1,2}(0)$. That is, paths
$\gamma _{1,1}$ and $\gamma _{1,2}$ are closed in
the corresponding complex planes ${\bf C}_1={\bf C}$
and ${\bf C}_2=\bf C$ embedded into $\bf H$.
In view of the conditions of the theorem, $\gamma _{1,1}$
for each $u$ and $\gamma _{1,2}$ for each $t$ corresponding to
$z\in \bf H$ are contained in subsets $\pi _{2,u}(U)$ and
$\pi _{1,t}(U)$ respectively which are simply connected.
Hence subsets $\Omega _{1,1}$ and $\Omega _{1,2}$ exist
in ${\bf C}_1$ and in ${\bf C}_2$ such that $\partial
\Omega _{1,1}=\gamma _{1,1}$ and $\partial \Omega _{1,2}=
\gamma _{1,2}$ and $\Omega _{1,1}\subset \pi _{2,u}(U)$
and $\Omega _{1,2}\subset \pi _{1,t}(U)$ for each $t$ and $u$
corresponding to $z\in U$ such that $\Omega _{1,1}$ and
$\Omega _{1,2}$ are simply connected in ${\bf C}_1$
and in ${\bf C}_2$ respectively.
It may easily be seen, taking into account \S 2.8, that
this integral can be considered as the integral of a real
differential $1$-form along the path $\gamma $ in $\bf R^4$.
To these integrals $Q$ and $\tilde Q$ the classical (generalized)
Stokes theorem can be applied (see Theorem V.1.1 \cite{weint}).
In view of the Hurewicz isomorphism theorem (see \S 7.5.4
\cite{span}) $H_q(U,x)=0$ for each $x\in U$ and each
$q<4$, hence $H^l(U,x)=0$ for each $l\ge 1$.
\par If $f: Y\to V$ is continuous, then $r\circ f: Y\to
\Omega $ is continuous, if $f$ is onto $V$, then
$r\circ f$ is onto $\Omega $, where $r: V\to \Omega $ is a retraction,
$V$, $Y$ and $\Omega $ are topological spaces.
The topological space $U$ is metrizable, hence for each
closed subset $\Omega $ in $U$ there exists a canonical closed subset
$V\subset U$ such that $V\supset \Omega $ and $\Omega $
is a retraction of $V$, that is, there exists a continuous mapping
$r: V\to \Omega $, $r(z)=z$ for each $z\in \Omega $
(see \cite{eng} and Theorem 7.1 \cite{isbell}).
Therefore, if $V$ is a $3$-connected canonical closed subset
of $U$ and $\Omega $ is a two dimensional $C^0$-manifold
such that $\Omega $ is a retraction of $V$, then $\Omega $
is simply connected, since each continuous mapping
$f: S^k\to \Omega $ with $k\le 1$ has a continuous extension
$f: {\bf R^{k+1}}\to V$ and $r\circ f: {\bf R^{k+1}}\to \Omega $
is also a continuous extension of $f$ from $S^k$
on $\bf R^{k+1}$.
\par From $3$-connectedness of $U$ it follows, that there are
two dimensional real differentiable manifolds $\Omega _j$ contained in
$U$ such that $\partial \Omega _j=\gamma $.
This may be lightly seen by considering partitions ${\sf Z}_n$ of $U$
by $S^n_{l,k}\cap U$ and taking $n\to \infty $, where
$S^n_{l,k}$ are parallelepipeds in $\bf H$ with ribs of length
$n^{-1}$, $l$, $k$ and $n\in \bf N$, two dimensional faces
$\mbox{ }_1S^n_l$ and $\mbox{ }_2S^n_k$ of $S^n_{l,k}=
\mbox{ }_1S^n_l\times \mbox{ }_2S^n_k$ are parallel to
$\bf C_1$ or $\bf C_2$ respectively such that there exists a sequence
of paths $\gamma _n$ converging to $\gamma $ relative to $|*|_{\bf H}$
and a sequence of (continuous) two dimensional $C^0$-manifolds
$\Omega ^n_j$ with $\partial \Omega ^n_j=\gamma ^n$,
$\Omega ^n_j\subset \bigcup_{l,k}[(\partial \mbox{ }_1S^n_l)\times (\partial
\mbox{ }_2S^n_k)]$.
Choose $\Omega _1$ and $\Omega _2$ orientable and of class $C^3$ as Riemann
manifolds such that taking their projections on $\bf C_1$ and
$\bf C_2$ the corresponding paths $\gamma _{1,1}$,
$\gamma _{1,2}$ and regions $\Omega ^j_{1,1}$ and $\Omega ^j_{1,2}$
in ${\bf C}_1$ and in ${\bf C}_2$ satisfy the conditions mentioned above
in this proof, where $j=1$ or $j=2$ for $\Omega _1$ or
$\Omega _2$ respectively. In this situation the abstract Stokes theorem
is applicable. In view of the Fubini theorem we obtain
$Q\tilde Q=\int_{\Omega _1\times \Omega _2} \eta (z_1,\tilde z_1)
\wedge \tilde \eta (z_2,\tilde z_2)$,
where $\eta =d((\partial g(z,\tilde z)/\partial z).dz)$
is the $2$-differential form, $z_1\in \Omega _1$,
$z_2\in \Omega _2$.
The function $g$ is in $C^{2,1}(U,{\bf H})$
and $p{\tilde p}={\tilde p}p$ for each $p\in \bf H$, hence
$(\partial ^2g(z,\tilde z)/\partial z\partial {\tilde z}).(h_1,h_2):=$
$(\partial [(\partial g(z,\tilde z)/\partial z).h_1]/\partial
{\tilde z}).h_2=$
$(\partial ^2g(z,\tilde z)/\partial {\tilde z}\partial z).(h_2,h_1)$
for each $h_1$ and $h_2$ in $\bf H$, in particular for $h_1=h_2=I$.
Due to Condition $(2.1)$ there is the equality
$(\partial f(z,\tilde z)/\partial {\tilde z})=0$, hence
$(\partial ^2g(z,\tilde z)/\partial z\partial {\tilde z})=$
$(\partial ^2g(z,\tilde z)/\partial {\tilde z}\partial z)=0$
and inevitably
$g=p(z)+q(\tilde z)$, where functions $p$ and $q$ are of class
$C^{2,1}$ such that $\partial p/\partial {\tilde z}=0$ and
$\partial q(\tilde z)/\partial z=0$.
This is evident in the class of polynomial functions,
that is dense in $C^{2,1}(W,{\bf H})$ for each compact canonical closed
set $W$ contained in $\bf H$ such that $W\subset U$. Hence it is true in
$C^{2,1}$ also. Therefore, $\partial _zg(z,\tilde z)=\partial _zp(z)$.
This means, that $\partial _zg(z,\tilde z)=dp(z)$, since
$\partial _{\tilde z}p=0$ (see Equation $(2.15)$).
Then $\eta =d^2(p)=0$, since $d=\mbox{ }_{\bf R}d$ and $d^2=0$.
Therefore, $Q=0$.
\par {\bf 2.12. Definitions.} A continuous function
on a domain $U$ in $\bf H$ such that $\emptyset \ne
Int (U)\subset U\subset cl(Int(U))$ and $\int_{\gamma }fdz=0$
for each rectifiable closed path $\gamma $ in $U$, then
$f$ is called quaternion integral holomorphic (on $U$).
\par If $f$ is a superdifferentiable function
on $U$ such that it satisfies Condition $(2.1)$, then
it is called quaternion holomorphic (on $U$).
\par Let $B(a,0,R,{\bf H})$ be a disk in $\bf H^n$,
then the completion of the space of all functions
having decomposition $(2.5.i,ii)$ with respect to
$z$ only, and with $n_0\ge 0$ relative to the norm
$\| * \|_{\omega }$ from \S 2.6, is denoted by $C^{\omega }_z
(B(a,0,R,{\bf H}),{\bf H})$. It is $\bf R$-linear.
Then $C^{\omega }_z(U,{\bf H})$
denotes the space of all continuous functions $f$ on $U$
with values in $\bf H$ such that for each $a\in U$
there are $R=R(f)>0$ and $g\in C^{\omega }_z(B(a,0,R,{\bf H}),{\bf H})$
with the restriction $g|_U=f$.
If $f\in C^{\omega }_z(U,{\bf H})$, then it is called
quaternion locally $z$-analytic (on $U$).
\par {\bf 2.13. Corollary.} {\it Let $f$ be a quaternion holomorphic
function on an open $3$-connected domain $U$ in $\bf H$
such that $\pi _{1,t}(U)$ and $\pi _{2,u}(U)$ are simply connected
in ${\bf C}_1$ for each $t$ and $u\in \bf C$ for which there exists $z\in U$,
then $f$ is quaternion integral holomorphic.}
\par This follows immediately from Theorem 2.11.
\par {\bf 2.14. Definition.} Let $U$ be a subset of $\bf H$
and $\gamma _0: [0,1]\to \bf H$ and $\gamma _1: [0,1]\to \bf H$
be two continuous paths. Then $\gamma _0$ and $\gamma _1$ are called
homotopic relative to $U$, if there exists a continuous mapping
$\gamma : [0,1]^2\to U$ such that $\gamma ([0,1],[0,1])\subset U$
and $\gamma (t,0)=\gamma _0(t)$ and $\gamma (t,1)=\gamma _1(t)$ for each
$t\in [0,1]$.
\par {\bf 2.15. Theorem.} {\it Let $W$ be an open subset in $\bf H$
and $f$ be a quaternion holomorphic function on $W$ with values in $\bf H$.
Suppose that there are
two rectifiable paths $\gamma _0$ and $\gamma _1$ in $W$ with common
initial and final points ($\gamma _0(0)=\gamma _1(0)$ and $\gamma _0(1)
=\gamma _1(1)$) homotopic relative to $U$, where $U$ is a $3$-connected
subset in $W$ such that $\pi _{1,t}(U)$ and $\pi _{2,u}(U)$
are simply connected in $\bf C$ for each $t$ and $u\in \bf C$
for which there exists $z\in \bf U$. 
Then $\int_{\gamma _0}fdz=\int_{\gamma _1}fdz$.}
\par {\bf Proof.} A homotopy of $\gamma _0$ with $\gamma _1$ realtive to $U$
implies homotopies of $(\gamma _0)_{1,j}$ with $(\gamma _1)_{1,j}$
relative to $\pi _{1,t}(U)$ and $\pi _{2,u}(U)$ in $\bf C$ for
$j=1$ and $j=2$ respectively for each $t$ and $u\in \bf C$ for which
there exists $z\in U$.
Consider a path $\zeta $ such that $\zeta (t)=
\gamma _0(2t)$ for each $0\le t\le 1/2$ and $\zeta (t)=\gamma _1(2-2t)$
for each $1/2 \le t\le 1$. Then $\zeta $ is a closed path
contained in a $U$. In view of Theorem 2.11 $\int_{\zeta }f(z)dz=0$.
On the other hand, $\int_{\zeta }f(z)dz=\int_{\gamma _0}f(z)dz-
\int_{\gamma _1}f(z)dz$, consequently, $\int_{\gamma _0}f(z)dz=
\int_{\gamma _1}f(z)dz$.
\par {\bf 2.16. Theorem.} {\it Let $f$ be a quaternion locally
$z$-analytic function on an open domain $U$ in $\bf H^n$, then $f$ is
quaternion holomorphic on $U$.}
\par {\bf Proof.} From the definition of the superdifferential
we get $(\partial z^n/\partial z).h=\sum_{k=0}^{n-1}z^khz^{n-k-1}$
and $\partial z^n/\partial {\tilde z}=0$.
Using the formula of the superdifferential for a product of functions,
from \S 2.7 we obtain, that each $f$ of the form
$(2.5.i,ii)$ is superdifferentiable
by $z$ when $n_0\ge 0$ and $m=0$ in $(2.5.ii)$ and hence satisfies
Condition $(2.1)$. Using the norm $\| *\|_{\omega }$-convergence
of series with respect to $z$ for a given $f\in C^{\omega }(U,{\bf H})$
we obtain for each $a\in U$, that there exists its neighbourhood $W$,
where $f$ is quaternion holomorphic, hence $f$ is quaternion
holomorphic on $U$.
\par {\bf 2.17. Note.} In the next section it is shown that
a quaternion holomorphic function is infinite differentiable;
furthermore, under suitable conditions equivalences between the
properties of quaternion holomorphicity, quaternion integral holomorphicity
and quaternion local $z$-analyticity, will be proved there too.
Integral $(2.6)$ may be generalized for a continuous function
$q: U\to \bf H$ such that $V(q\circ \gamma )<\infty $.
Substituting $\Delta z_k$ on $q(z_{k+1})-q(z_k)=:\Delta q_k$
in $(2.7)$ we get
$$(2.6')\quad \int_{\gamma }f(z,{\tilde z})dq(z):=
\lim_PI(f,q\circ \gamma ;P),
\mbox{ where}$$
$$(2.7')\quad I(f,q\circ \gamma ;P)=\sum_{k=0}^{q-1}{\hat f}(z_{k+1},
{\tilde z}_{k+1}).(\Delta q_k).$$
In paticular, if $\gamma \in C^1$ and $q$ is quaternion holomorphic
on $U$, also $f(z,{\tilde z})=(\partial g/\partial z).I$, where
$g\in C^{1,0}(U,{\bf H})$, then
$$\int_{\gamma }f(z,{\tilde z})dq(z)=\int_0^1
(\partial g/\partial z).((D_zq(z)|_{z=\gamma (s)}).\gamma '(s))ds$$
and $V(\gamma )\le \int_0^1 |\gamma '(s)| ds$.
\par Let $f: U\to \bf H$, where $U$ is an open subset of $\bf H^n$.
If there exists a quaternion holomorphic function
$g: U\to \bf H$ such that $g'(z).I=f(z)$
for each $z\in U$, then $g$ is called a primitive of $f$.
\par {\bf 2.18. Proposition.} {\it Let $U$ be an open connected subset
of $\bf H^n$ and $g$ be a primitive of $f$ on $U$,
then a set of all primitives of $f$ is:
$\{ h: h=g+C, C=const \in {\bf H} \} $.}
\par {\bf Proof.} Suppose $h'(z)=0$ for each $z\in U$, then
consider $q(s):=h((1-s)a+sz)$ for each $s\in [0,r]$, where
$a$ is a marked point in $U$ and $B(a,r,{\bf H})$ is a ball
contained in $U$, $r>0$, $z\in B(a,r,{\bf H})$. Then
$q$ is correctly defined and $q(0)=q(1)$. Therefore,
the set $V:=\{ z\in U: h(z)=h(a) \} $ is open in $U$, since
with each point $a$ it contains its neighbourhood.
On the other hand, it is closed due continuity of $h$,
hence $V=U$, since $U$ is connected, consequently, $h=const $ on $U$.
\section{Meromorphic functions and their residues.}
At first we define and describe the exponential and
the logarithmic functions of quaternion variables
and then apply them to the investigation of quaternionic residues.
\par {\bf 3.1. Note and Definition.} Let $z\in \bf H$,
then
$$(3.1.)\quad \exp (z):=\sum_{n=0}^{\infty }z^n/n! .$$
This definition is correct, since real numbers commute with
quaternions. If $|z|\le R<\infty $, then the series $(3.1)$
converges, since
$|exp (z)| \le \sum_{n=0}^{\infty } |z^n/n!|\le \exp (R)<\infty .$
Therefore, $\exp $ is the function defined on $\bf H$
with values in $\bf H$. The restriction of $\exp $ on
the subset ${\bf Q_d}:=\{ z: z\in {\bf H}, z={{t\quad 0}\choose
{0\quad \bar t}},$ $t\in {\bf C} \}$ is commutative, but in general
two quaternions $z_1$ and $z_2$ do not commute and
on ${\bf H}^2$ the function $\exp (z_1+z_2)$ does not
coincide with $\exp (z_1)\exp (z_2)$.
\par {\bf 3.2. Proposition.} {\it Let $z\in \bf H$ be written
in the form $z=vI+s(wJ+xK+yL)$ with real $v$, $w$, $x$, $y$ and $s$
with $w^2+x^2+y^2\ne 0$,
then
$$(3.2) \quad \exp (z)=\exp (v) \{ \cos (s\phi )I+ i(\sin (s\phi )/
\phi )
{{w\quad (y-ix)}\choose {(y+ix)\quad -w}} \} ,$$
where $\phi :=(w^2+x^2+y^2)^{1/2}$.}
\par {\bf Proof.} Since the unit matrix $I$ commutes with
$J$, $K$ and $L$, then $\exp (vI+z)=\exp (v)\exp (z)$ for each
$v\in \bf R$ and each $z\in \bf H$. Consider $z=s(wJ+xK+yL)$
with real $s$, $w$, $x$ and $y$ such that $w^2+x^2+y^2=1$,
then
$$z=is{{w\quad (y-ix)}\choose {(y+ix)\quad -w}},$$ where
$i=(-1)^{1/2}$. Denote $y-ix=:m$, then
$${{w\quad m}\choose {\bar m\quad -w}}^{2k}=(w^2+|m|^2)^kI$$
for each $k\in \bf N$. On the other hand $|m|^2=x^2+y^2$,
hence $w^2+|m|^2=1$, consequently,
$$\sum_{k=0}^{\infty }z^k/k!=\sum_{k=0}^{\infty }
(-1)^ks^{2k}I/(2k)!
+i\sum_{k=0}^{\infty }(-1)^ks^{2k+1}{{w\quad m}\choose {\bar m\quad -w}}
/(2k+1)!$$
$$=\cos (s)I+i\sin (s) {{w\quad m}\choose {\bar m\quad -w}}.$$
The particular case $s=0$ corresponds to $exp(0)=1$.
From this Formula $(3.2)$ follows.
\par {\bf 3.3. Corollary.} {\it If $z\in \bf H$
is written in the form $z=vI+wJ+xK+yL$ with real $v$, $w$, $x$ and $y$,
then $|\exp (z)|=\exp (v)$.}
\par {\bf Proof.} If $w^2+x^2+y^2=0$ this is evident.
Suppose $w^2+x^2+y^2\ne 0$. In view of Formula $(3.2)$
$$(3.3)\quad exp(z)=exp(v)A\mbox{, where }A={{(\cos (\phi )+i\sin (\phi )w)
\quad  \sin (\phi ) (x+iy)}\choose
{\sin (\phi )(-x+iy)\quad (\cos (\phi ) -i\sin (\phi )w)}}.$$
Since $A\in \bf H$, then $|A|^2=det (A)=1$
and inevitably $|\exp (z)|=\exp (v)$.
\par {\bf 3.4. Corollary.} {\it The function $\exp (z)$
on the set ${\bf H_i}:=\{ z: z\in {\bf H}, Re (z_{1,1})=0 \} $ is periodic
with three generators of periods $J$, $K$ and $L$
such that $\exp (z(1+2\pi n/|z|))=\exp (z)$ for each
$0\ne z\in \bf H_i$ and each integer number $n$.
If $z\in \bf H$ is written in the form $z=2\pi sM$, where $M=wJ+xK+yL$,
with real $w$, $x$ and $y$ such that $w^2+x^2+y^2=1$,
then $\exp (z)=1$ if and only if $s\in {\bf Z}$.}
\par {\bf Proof.} In view of Formula $(3.2)$
$\exp (z)=1$ for a given $z\in \bf H_i$ if and only if
$\cos (s\phi )=1$ and $\sin (s\phi )=0$, that is equivalent to
$s\in \{ 2\pi n: n\in {\bf Z} \} $, since $\phi =1$
by the hypothesis of this corollary.
The particular cases of Formula $(3.2)$
are either $w\ne 0$, $x=y=0$; or $w=y=0$ and $x\ne 0$;
or $w=x=0$ and $y\ne 0$, hence $J$, $K$ and $L$ are the three
generators for the periods of $\exp $.
\par {\bf 3.5. Corollary.} {\it The function $\exp $
is the epimorphism from $\bf H_i$ on the three-dimensional
quaternion unit sphere $S^3(0,1,{\bf H}):=\{ z: z\in {\bf H},
|z|=1 \} $. }
\par {\bf Proof.} In view of Corollary 3.3 the image
$\exp ({\bf H_i})$ is contained in $S^3(0,1,{\bf H})$.
The sphere $S^3(0,1,{\bf H})$ is characterized by the condition
$v_1^2+w_1^2+x_1^2+y_1^2=1$. In view of Formula $(3.2)$
we have $v_1=\cos (s)$, $w_1=\sin (s)w$, $x_1=
\sin (s)x$ and $y_1=\sin (s)y$, where $s\in \bf R$ and
$w^2+x^2+y^2=1$. Vice versa let $z_1\in S^3(0,1,{\bf H})$.
For each $v_1\in [-1,1]$ there exists
$s=arcsin (v_1)$ such that $v_1=\cos (s)$ and $w_1^2+x_1^2+y_1^2=\sin ^2(s)$.
The case $\sin (s)=0$ corresponds to $v_1=1$ and others coordinates
equal to zero, hence $z_1=\exp (0)$. If $\sin (s)\ne 0$
there are $w=w_1/\sin (s)$, $x=x_1/\sin (s)$ and $y=y_1/\sin (s)$,
consequently, $\exp (z)=z_1$ in this case too.
Therefore, $\exp $ is an epimorphism of $\bf H_i$
on $S^3(0,1,{\bf H})$. 
\par {\bf 3.6. Corollary.} {\it Each quaternion has a polar decomposition
$$(3.4)\quad z=\rho \exp (2\pi (\phi _1J+\phi _2K+\phi _3L)),$$
where $\phi _j\in [-1,1]$ for each $j=1,2,3$,
$\phi _1^2+\phi _2^2+\phi _3^2=1$, $\rho :=|z|.$}
\par {\bf Proof.} This follows from Formula $(3.2)$ and
Corollary 3.5.
\par {\bf 3.7. Note.} In the noncommutative quaternion case
there is the following relation for $\exp $ and its (right) derivative:
$$(3.5) \quad \exp (z)'.h=\sum_{n=1}^{\infty }\sum_{k=0}^{n-1}
z^khz^{n-k-1}/n!,$$
where $z$ and $h\in \bf H$. In particular,
$$(3.6)\quad \exp (z)'.vI=v \exp (z)$$ for each $v\in \bf R$,
but generally not for all $h\in \bf H$.
In view of \S 2.6 the derivative $(3.5)$ reduces to the form
of the definition of superdifferentiability given in \S 2.1.
The function $\exp $ is periodic on $\bf H$, hence
the inverse function denoted by $Ln$ is defined only locally.
Consider the space $\bf R^3$ of all variables $w$, $x$ and $y$
for which $\exp $ is periodic on $\bf H$. The condition
$w^2+x^2+y^2=1$ defines in $\bf R^3$ the unit sphere $S^2$.
The latter has a central symmetry element $C$ for the transformation
$C(w,x,y)=(-w,-x,-y)$. Consider a subset $P=\bigcup_{k=1}^4P_k$ of $S^2$
of all points characterized by the conditions:
$P_1:=\{ (w,x,y)\in S^2: w\le 0, x\le 0, y\le 0 \} ,$
$P_2:=\{ (w,x,y)\in S^2: w\ge 0, x\le 0, y\le 0 \} ,$
$P_3:=\{ (w,x,y)\in S^2: w\le 0, x\ge 0, y\le 0 \} ,$
$P_4:=\{ (w,x,y)\in S^2: w\le 0, x\le 0, y\ge 0 \} ,$
then $P\cup CP=S^2$ and the intersection $P\cap CP$
is one dimensional over $\bf R$.
This sphere $S^2$ corresponds to
the embedding $\theta _1: (w,x,y)\hookrightarrow (0,w,x,y)\in \bf R^4$.
Consider the embedding of $\bf R^4$ into $\bf H$
given by $\theta _2: (v,w,x,y)\hookrightarrow vI+wJ+xK+yL\in \bf H$.
This yields the embedding $\theta :=\theta _2\circ \theta _1$
of $S^2$ in $\bf H$.
Each unit circle with the center $0$ in $\bf H$ intersects
the equator $\theta (S^2)$ of the unit sphere $S^3(0,1,{\bf H})$.
Join each point $(wJ+xK+yL)$ on $\theta (S^2)$ with the zero point in
$\bf H$ by a line $\{ s(wJ+xK+yL): s\in {\bar {\bf R}}_+ \} $,
where ${\bar {\bf R}}_+:=\{ s\in {\bf R}: s\ge 0 \} $.
This line crosses a circle embedded into $S^3(0,1,{\bf H})$,
which is a trace of a circle $\{ \exp (2\pi s(wJ+xK+yL)): s\in [0,1] \} $
of radius $1$
in $\bf H$. Therefore, $\psi (s):= \exp (vI+2\pi s(wJ+xK+yL))$
as a function of $(v,s)$ for fixed $(w,x,y)\in S^2$ defines a bijection
of the domain $X\setminus \{ s(wJ+xK+yL): s\in {\bar {\bf R}}_+ \} $
onto its image,
where $X$ is $\bf R^2$ embedded as $(v,s)\hookrightarrow
(vI+s(wJ+xK+yL)) \in \bf H$. This means, that $Ln (z)$ is correctly
defined on each subset $X\setminus \{ s(wJ+xK+yL): s\in {\bar {\bf R}}_+ \} $
in $\bf H$. The union $\bigcup_{(w,x,y)\in P}
\{ s(wJ+xK+yL): s\in {\bar {\bf R}}_+ \} $ produces the three dimensional
(over $\bf R$) subset $Q:=\bigcup_{k=1}^4Q_k$, where
$Q_k:=\theta (S_k)$, $S_1:=\{ (w,x,y)\in {\bf R^3} \} :
w\le 0, x\le 0, y\le 0 \} $,
$S_2:=\{ (w,x,y)\in {\bf R^3}: w\ge 0, x\le 0, y\le 0 \} $,
$S_3:=\{ (w,x,y)\in {\bf R^3}: w\le 0, x\ge 0, y\le 0 \} $,
$S_4:=\{ (w,x,y)\in {\bf R^3}: w\le 0, x\le 0, y\ge 0 \} $,
${\bar {\bf R}}_+:=[0,\infty )$.
Then, on the domain ${\bf H}\setminus Q$, the function $\exp (z)$
defines a bijection with image $\exp ({\bf H}\setminus Q)$
and its inverse function $Ln (z)$ is correctly defined
on ${\bf H}\setminus \exp (Q)$.
By rotating ${\bf H}\setminus Q$ one may produce other
domains on which $Ln$ can be defined as
the univalued function (that is, $Ln(z)$ is one point in $\bf H$),
but not on the entire $\bf H$. This means that $Ln(z)$ is a
locally bijective function.
We have elementary identities $\cos (2\pi -\phi )=\cos (\phi )$
and $\sin (2\pi -\phi )=-\sin (\phi )$ for each $\phi \in \bf R$.
If $0<\phi <2\pi $, then $w_1\sin (\phi )/\phi  =w_2\sin (2\pi -\phi )/
(2\pi -\phi )$ if and only if $w_1=-\phi w_2/(2\pi -\phi ).$
To exclude this ambiguity we put in Formula $(3.2)$
$\phi \ge 0$ such that $\phi =(w^2+x^2+y^2)^{1/2}$ is
the nonegative (arithmetical) branch of the square root function
on ${\bar {\bf R}}_+$ and $w\ge 0$.
Therefore, $Ln(\exp (z))=z$ on ${\bf H}\setminus Q$, hence
using Formulas $(3.3,3.4)$ we obtain the multivalued function
$$(3.7)\quad Ln (z)=ln (|z|)+Arg(z), \mbox{ where }
Arg(z):=arg (z)+2\pi sM$$
on ${\bf H}\setminus \{ 0 \} $, where $ln$ is the usual
real logarithm on $(0,\infty )$, $s\in \bf Z$,
$$|z|\exp (2\pi arg (z))=z,\quad arg (z):=w_zJ+x_zK+y_zL,\quad
(w_z,x_z,y_z)\in {\bf R^3},$$
$w_z^2+x_z^2+y_z^2<1$, $w_z\ge 0$,
$M=wJ+xK+yL$ is any unit vector (that is, $|M|=1$) in $\bf H$ commuting with
$arg (z)\in \bf H$, $arg (z)$ is uniquely defined by such
restriction on $(w_z,x_z,y_z)$, for example, $M=\zeta arg (z)$
for any $\zeta \in \bf R$, when $arg (z)\ne 0$.
\par For each fixed
$M=wJ+xK+yL$ $\quad \exp (sM)$ is a one-parameter family of
special unitary transformations (that is, $det (\exp (sM))=1$)
of $\bf H$ (that induces rotations of
the Euclidean space $\bf R^4$), that is, $\exp (sM)\eta \in \bf H$
for each $\eta \in \bf H$, where $\bf H$ as the linear space
over $\bf R$ is isomorphic with $\bf R^4$. On the other hand, there
are special unitary transformations of $\bf H$ for which $s=\pi /2+\pi k$,
but $M$ is variable with $|M|=1$, where
$k\in \bf Z$, then $\exp (z)=(-1)^k A,$ where
$A={{iw\quad (x+iy)}\choose {(-x+iy)\quad -iw}}$ (see Formula $(3.3)$).
To each closed curve $\gamma $ in $\bf H$ there corresponds
a closed curve $P_{\xi }(\gamma )$ in a $\bf R$-linear subspace
$\xi \ni 0$, where $P_{\xi }$ is a projection on $\xi $, for example,
$P_{{\bf R}I\oplus {\bf R}J}(z)=(z-JzJ)/2=vI+wJ$
for $\xi ={\bf R}I\oplus {\bf R}J$,
$P_{{\bf R}I\oplus {\bf R}K}(z)=(z-KzK)/2=vI+xK$,
$P_{{\bf R}I\oplus {\bf R}L}(z)=(z-LzL)/2=vI+yL$,
$P_{{\bf R}J\oplus {\bf R}K}(z)=(z+LzL)/2=wJ+xK$,
$P_{{\bf R}J\oplus {\bf R}L}(z)=(z+KzK)/2=wJ+yL$,
$P_{{\bf R}K\oplus {\bf R}L}(z)=(z+JzJ)/2=xK+yL$,
$P_{{\bf R}J\oplus {\bf R}K\oplus {\bf R}L}(z)=(3z+JzJ+KzK+LzL)/4=
wJ+xK+yL$ for $\xi ={\bf R}J\oplus {\bf R}K\oplus {\bf R}L$, etc.
Particular cases of special unitary transformations also correspond to
$w=0$ or $x=0$ or $y=0$ for $M\ne 0$. To each closed curve
$\gamma $ in $\bf H$ and each quaternions $a$ and $b$ with
$ab\ne 0$ there corresponds a closed curve $a\gamma b$ in $\bf H$,
for example, for $a=J$ and $b=K$ there is the identity
$JzK=vL-wK-xJ+yI$ for each $z=vI+wJ+xK+yL$ in $\bf H$.
\par Instead of the Riemann two dimensional surface of the
complex logarithm function we get the four dimensional
manifold $W$, that is, a subset of $Y^{\aleph _0}:=
\prod_{i\in \bf Z}Y_i,$ where $Y_i= Y$ for each $i$, such that
each $Y$ is a copy of ${\bf H}$ embedded into ${\bf H}\times \bf R^3$
and cut by a three dimensional submanifold
$Q$ and with diffeomorphic bending of a neighbourhood of $Q$
such that two three dimensional edges $Q_1$ and $Q_2$ of $Y$ diffeomorphic
to $Q$ do not intersect outside zero, $Q_1\cap Q_2= \{ 0 \} ,$
that is, the boundary $\partial Q$ is also cut everywhere outside
zero. We have $\partial Q= \partial Q^w\cup \partial Q^x\cup \partial Q^y$,
where $\partial Q^w:=\{ \theta (w,x,y): w=0, (w,x,y)\in \bigcup_{k=1}^4S_k \} $,
$\partial Q^x:=\{ \theta (w,x,y): x=0, (w,x,y)\in \bigcup_{k=1}^4S_k \} $
and
$\partial Q^y:=\{ \theta (w,x,y): y=0, (w,x,y)\in \bigcup_{k=1}^4S_k \} $.
This means, that $\partial Q^w=\{ z=xK+yL: (x,y)\in {\bf R^2},$
$x\mbox{ and } y $ $\mbox{ are not simultaneously positive} \} $.
Similarly for $\partial Q^x$ and $\partial Q^y$
with $z=wJ+yL$ and $z=wJ+xK$
instead of $z=xK+yL$ respectively.
To exclude rotations in each subspace $vI+s(aK+bL)$
isomorphic with $\bf R^2$ and embedded into ${\bf R}I+\partial Q^w$
and similarly for $vI+s(aJ+bL)$ and $vI+s(aJ+bK)$
we have cut $\partial Q$, where $v$, $s\in \bf R$ are variables
and $a,$ $b$ are two real constants such that $ab\le 0$, $a^2+b^2>0$.
Then in ${\bf H}\times \bf R^3$ two copies
$Y_i$ and $Y_{i+1}$ are glued by the equivalence relation
of $Q_{2,i}$ with $Q_{1,i}$ via the segments $\{ s_{l,i}(wJ+xK+yL): s_{l,i}
\in {\bar {\bf R}}_+  \} $
such that $s_{1,i+1}=s_{2,i}$ for each $s_{l,i}\in {\bar {\bf R}}_+$
and each given real $(w,x,y)\in P$ with $w^2+x^2+y^2=1$.
This defines the four dimensional manifold $W$ embedded into ${\bf H}
\times \bf R^3$ and $Ln: {\bf H}\setminus \{ 0 \} \to W$ is the
univalued function, that is, $Ln(z)$ is a singleton
in $W$ for each $z\in {\bf H}\setminus \{ 0 \}$.
\par {\bf 3.8. Theorem.} {\it The function $Ln$ is quaternion
holomorphic on any domain $U$ in $\bf H$ obtained by a
quaternion holomorphic diffeomorphism of ${\bf H}\setminus Q$
onto $U$.
Each path $\gamma $ in $\bf H$ such that $\gamma (s)=r\exp (2\pi
sn(wJ+xK+yL))$ with $s\in [0,1]$, $n\in {\bar {\bf R}}_+$,
$w^2+x^2+y^2=1$ is closed in $\bf H$ if and only if $n\in \bf N$,
where $r>0$. In this case
$$(3.8)\quad \int_{\gamma }z^{-1}dz=
\int_{\gamma }d(Ln z)=2\pi n(wJ+xK+yL).$$}
\par {\bf Proof.} If $U$ and $V$ are two open subsets in
$\bf H$ and $g: V\to U$ is a quaternion holomorphic
diffeomorphism of $V$ onto $U$
and $f$ is a quaternion holomorphic function on $V$,
then $f\circ g^{-1}$ is quaternion holomorphic function
on $U$, since $(f\circ g^{-1})'(z).h=(f'(\zeta )|_{\zeta =
g^{-1}(z)}.(g^{-1}(z))'.h$ for each $z\in U$ and each
$h\in \bf H$ and $\partial (f\circ g^{-1}(z))/\partial {\tilde z}=
(\partial f(\zeta )/\partial {\tilde \zeta })|_{\zeta =g^{-1}(z)}.
(\partial {\tilde g}^{-1}(z)/\partial {\tilde z})+$
$(\partial f(\zeta )/\partial \zeta )|_{\zeta =g^{-1}(z)}.
(\partial g^{-1}(z)/\partial {\tilde z})=0$.
Since $\exp $ is the diffeomorphism
from ${\bf H}\setminus Q$ onto ${\bf H}\setminus \exp (Q)$, we have that
$Ln $ is quaternion holomorphic on  ${\bf H}\setminus Q$
and on each of its quaternion holomorphic images after choosing
a definite branch of the multivalued function $Ln (z)$ (see Formula $(3.7)$).
\par A path $\gamma $ is defined for each $s\in \bf R$
not only for $s\in [0,1]$ due to the existence of $\exp $.
In view of Formula $(3.2)$ a path $\gamma $ is closed
(that is, $\gamma (s_0)=\gamma (s_0+1)$ for each $s_0\in \bf R$)
if and only if $\cos (2\pi n)=\cos (0)=1$ and $\sin (2\pi n)=0$, that is,
$n\in \bf N$.
\par From the definition of the line integral
we get the equality: $\int_{\gamma }d(Lnz)=
\int _0^1(Ln z)'.\gamma '(s)ds$. Considering integral sums
by partitions $P$ of $[0,1]$ and taking the limit by the family
of all $P$ we get, that $\int_{\gamma }d(Ln z)=Arg(\gamma (1))-
Arg (\gamma (0))$ for a chosen branch of the function $Arg (z)$
(see Formula $(3.7)$). Therefore, $\int_{\gamma }d(Ln z)=2\pi n(wJ+xK+yL).$
\par Since $z$ commutes with itself, we have:
$\exp (z)'.z=\exp (z)z.$ Therefore,
$\exp (Ln(z))'.I=(\partial \exp (\eta )/\partial \eta )|_{\eta =Ln(z)}
.(Ln(z))'.I=\exp (Ln(z))(Ln(z))'.I$, consequently,
$(Ln (z))'.I=\exp (-Ln(z))=z^{-1}$ and inevitably
$$\lim_PI(z^{-1},\gamma ;P)=\lim_P\sum_l{\hat z}_l^{-1}\Delta z_l
=\lim_P\Delta Ln(z_l)=\int_{\gamma }dLn(z),$$
hence $\int_{\gamma }z^{-1}dz=\int_{\gamma }dLn(z)$.
That is, $\int_{\gamma }d Ln (z)$ can be considered
as the definition of $\int_{\gamma }z^{-1}dz$.
\par {\bf 3.9. Theorem.} {\it Let $f$ be a continuous quaternion
holomorphic function on an open domain $U$ in $\bf H$.
If $(\gamma +z_0)$ and $\psi $ are presented as piecewise unions of paths
$\gamma _j+z_0$ and $\psi _j$ with respect to parameter
$s\in [a_j,b_j]$ and $s\in [c_j,d_j]$
respectively with $a_j<b_j$ and $c_j<d_j$ for each $j=1,...,n$
and $\bigcup_j[a_j,b_j]=\bigcup_j[c_j,d_j]=[0,1]$ 
homotopic relative to $U_j\setminus \{ z_0 \} $, where
$U_j\setminus \{ z_0 \} $ is a $3$-connected open domain in $\bf H$
such that $\pi _{1,t}(U_j\setminus \{ z_0 \})$ and
$\pi _{2,u}(U_j\setminus \{ z_0 \})$ 
are simply connected in $\bf C$ for each $t$ and $u\in \bf C$
for which there exists $z\in \bf H$ for each $j$.
If $(\gamma +z_0)$ and $\psi $ are
closed rectifiable paths in $U$ such that
$\gamma (s)=r\exp (2\pi s(wJ+xK+yL))$ with $s\in [0,1]$
and $w^2+x^2+y^2=1$ and $z_0\notin \psi $. Then
$$(3.9)\quad f(z)=(2\pi )^{-1}(\int_{\psi }f(\zeta )(\zeta -z)^{-1}d\zeta )
(wJ+xK+yL)^{-1}$$
for each $z\in U$ such that $|z-z_0|< \inf_{\zeta \in \psi ([0,1])}
|\zeta -z_0|$.}
\par {\bf Proof.} Join $\gamma $ and $\psi $ by a rectifiable path
$\omega $ such that $z_0\notin \omega $,
which is going in one direction and the opposite direction,
denoted $\omega ^-$, 
such that $\omega _j\cup \psi _j\cup \gamma _j\cup \omega _{j+1}$
is homotopic to a point relative to $U_j\setminus \{ z_0 \}$
for suitable $\omega _j$ and $\omega _{j+1}$, where $\omega _j$ joins
$\gamma (a_j)$ with $\psi (c_j)$ and $\omega _{j+1}$ joins
$\psi (d_j)$ with $\gamma (b_j)$ such that $z$ and $z_0\notin
\omega _j$ for each $j$. Then
$\int_{\omega _j} f(\zeta )(\zeta -z)^{-1}d\zeta = -
\int_{\omega _j^-}f(\zeta )(\zeta -z)^{-1}d\zeta $
for each $j$. In view of Theorem 2.15 there is the equality
$-\int_{\gamma ^-+z}f(\zeta )(\zeta -z)^{-1}d\zeta =
\int_{\psi }f(\zeta )(\zeta -z)^{-1}d\zeta .$
Since $\gamma +z$ is a circle around $z$ its radius $r>0$ can be
chosen so small, that $f(\zeta )=f(z)+ \alpha (\zeta ,z)$,
where $\alpha $ is a continuous function on $U^2$ such that
$\lim_{\zeta \to z}\alpha (\zeta ,z)=0$, then
$\int_{\gamma +z}f(\zeta )(\zeta -z)^{-1}d\zeta =$
$\int_{\gamma +z}f(z)(\zeta -z)^{-1}d\zeta +\delta (r)=$
$2\pi f(z) (wJ+xK+yL)+\delta (r),$
where $|\delta (r)|\le |\int_{\gamma }\alpha (\zeta ,z)
(\zeta -z)^{-1}d\zeta |$ $\le 2\pi \sup_{\zeta \in \gamma }
|\alpha (\zeta ,z)| C_1 \exp (C_2 r^6)$,
where $C_1$ and $C_2$ are positive constants (see Inequality
$(2.12)$), hence there exists $\lim_{r\to 0,
r>0}\delta (r)=0$. Taking the limit while $r>0$ tends to zero
yields the conclusion of this theorem.
\par {\bf 3.9.1. Corollary.} {\it Let $f$, $U$, $\psi $,
$z$ and $z_0$ be as in Theorem 3.9, then
$$|f(z)|\le \sup_{(\zeta
\in \psi , h\in {\bf H}, |h|\le 1)} |{\hat f} (\zeta ).h|.$$}
\par {\bf 3.10. Theorem.} {\it Let $f$ be a continuous function
on an open subset $U$ of $\bf H$.
If $f$ is quaternion integral holomorphic on $U$, then
$f$ is quaternion locally $z$-analytic on $U$.}
\par {\bf Proof.} Let $z_0\in U$ be a marked point and let
$\Gamma $ denotes the family of all rectifiable paths
$\gamma : [0,1]\to U$ such that $\gamma (0)=z_0$, then
$U_0=\{ \gamma (1): \gamma \in \Gamma \} $ is a connected component
of $z_0$ in $U$. Therefore, $g= \{ \gamma (1), \int_{\gamma }
f(z)dz \} $ is the function with the domain $U_0$.
Let $X$ be a compact metric space and $F$ be a function
continuous on $U\times X$ with values in $\bf H$ and for each
$p\in X$ let $f_p(z):=F(z,p)$ be quaternion holomorphic on $U$ by
$z\in U$. Define $G$ on $U^2\times X$ by $G(z,w,p):=[F(z,p)-F(w,p)]
(w-z)^{-1}$, $w\ne z$. Then $G(z,z,p)=(\partial f_p(z)/\partial z).I$.
It can be seen with the help of Formula $(3.9)$
that $G$ is continuous on $U^2\times X$,
since $G(b,c,q)-G(a,a,p)=\int_{\gamma }[(\partial f_q(z)/\partial z)-
(\partial f_p(a)/\partial a)].dz.(c-b)^{-1}$,
where $\gamma $ is a rectifiable curve such that
$\gamma (0)=b$, $\gamma (1)=c$. Moreover, $G$ is uniformly continuous
on $V^2\times X$ for each compact canonical closed subset $V$ in $\bf H$
such that $V\subset U$.
As in \S 2.15 it can be proved, that
$F(z):=\int_{\gamma }f(z)dz$, for each rectifiable $\gamma $ in
$U$, depends only on initial and final points.
This integral is finite, since $\gamma ([0,1])$ is contained
in a compact canonical closed subset $W\subset U$
on which $f$ is bounded.
Therefore, $(\partial \int_{z_0}^zf(\zeta )d\zeta /\partial z).h=
{\hat f}(z).h$ for each $z\in U$ and $h\in \bf H$,
$(\partial \int_{z_0}^zf(\zeta )d\zeta /\partial {\tilde z})=0$
for each $z\in U$ and $h\in \bf H$, where
$z_0$ is a marked point in $U$ such that $z$ and $z_0$ are in one
connected component of $U$.
In particular, ${\hat f}(z).I=f(z)$ for each $z\in U$.
Here $\hat f$ is correctly defined for each $f\in C^{1,0}(U,{\bf H})$
by continuity of the differentiable integral functional
on $C^0(U,{\bf H})$. In particular, ${\hat f}(z).I=f(z)$
for each $z\in U$. For a given $z\in \bf U$ choose
a neighbourhood $W$ satisfying the conditions of Theorem 3.9.
Then there exists a rectifiable path $\psi \subset W$ such that $F(z)$
is presented by Formula $(3.9)$. The latter integral is
infinite differentiable by $z$ such that
$$(3.10)\quad (\partial ^kF(z)/\partial z^k)=
k!(2\pi )^{-1}(\int_{\psi }f(\zeta )(\zeta -z)^{-k-1}d\zeta )
(w_0J+x_0K+y_0L)^{-1},$$
where $w_0$, $x_0$ and $y_0\in \bf R$ are fixed and $w_0^2+x_0^2+y_0^2=1$. 
In particular, we may choose a ball $W=B(a,R,{\bf H}):=
\{ \xi \in {\bf H}: |\xi -a|\le R \} \subset U$
for a sufficiently small $R>0$ and $\psi =\gamma +a$, where
$\gamma (s)=r\exp (2\pi s(w_0J+x_0K+y_0L))$ with $s\in [0,1]$, $0<r<R$.
If we prove that $F(z)$ is quaternion locally $z$-analytic, then
evidently its $z$-derivative $f(z)$ will also be  
quaternion locally $z$-analytic.
Consider $z\in B(a,r',{\bf H})$  with $0<r'<r$, then $|z-a|<|\zeta -a|$
for each $\zeta \in \psi $ and $(\zeta -a -(z-a))^{-1}=
(1-(\zeta -a)^{-1}(z-a))^{-1}(\zeta -a)^{-1}=$ $\sum_{k=0}^{\infty }
((\zeta -a)^{-1}(z-a))^k(\zeta -a)^{-1}$,
where $0\notin \psi $.
Therefore,
$$(3.11)\quad F(z)= (2\pi )^{-1}\sum_{k=0}^{\infty }\phi _k(z),$$
$$\mbox{where }\phi _k(z):=
(\int_{\psi }f(\zeta )((\zeta -a)^{-1}(z-a))^k(\zeta -a)^{-1}d\zeta )
(w_0J+x_0K+y_0L)^{-1}.$$
Thus $|\phi _k(z)|\le \sup_{\zeta \in \psi }|f(\zeta )|(r'/r)^{-k}$
for each $z\in B(a,r',{\bf H})$ and series $(3.11)$ converges
uniformly on $B(a,r',{\bf H})$. Each function $\phi _k(z)$
is evidently quaternion locally $z$-analytic on $B(a,r',{\bf H})$,
hence $F(z)$ is such too. Since for each $a\in U$
there is an $r'>0$, for which the foregoing holds,
it follows that $f(z)$ is the quaternion locally $z$-analytic function.
\par {\bf 3.11. Note.} Theorems 2.11, 2.15, 2.16, 3.10 and Corollary 2.13
establish the equivalence of notions of quaternion holomorphic,
quaternion integral holomorphic
and quaternion locally $z$-analytic classes of functions on domains
satisfying definite conditions.
Before, the notion of quaternion holomorphicity was defined relative to
a right superdifferentiation, similarly it can be defined relative
to a left superdifferentiation. Quaternion local $z$-analyticity shows,
that a function is quaternion holomorphic relative
to a right superdifferentiation if and only if it is
quaternion holomorphic relative to a left superdifferentiation.
\par In particular, if $f\in \mbox{ }_lC^{\omega }(U,{\bf H})$,
then evidently $F(z):=\int_{z_0}^zf(\zeta )d\zeta $ and
$(\partial f(\zeta )/\partial \zeta ).I$ belong to
$\mbox{ }_lC^{\omega }(U,{\bf H})$, where $z$ and $z_0\in U_0$,
$\zeta \in U$, $U_0$ is a connected component of $U$ open in
$\bf H$, since $(b_n{\hat \zeta }^n).\Delta \zeta =b_n(\partial
\zeta ^{n+1}/\partial \zeta ).\Delta \zeta $ for each $\zeta
\in \bf H$, $\Delta \zeta \in \bf H$, $n\in \bf N$, $b_n\in \bf H$.
\par {\bf 3.11.1. Definitions.} Let $U$ be an open subset in $\bf H$ and
$f\in C^0(U,{\bf H})$, then we say that $f$ possesses a primitive
$g\in C^1(U,{\bf H})$ if $g'(z).I=f(z)$ for each $z\in U$.
A region $U$ in $\bf H$
is said to be quaternion holomorphically simply connected
if every function quaternion holomorphic on it possesses
a primitive.
\par From \S 3.10 we get.
\par {\bf 3.11.2 Theorem.} {\it If $f\in C^{\omega }(U,{\bf H})$,
where $U$ is $3$-connected; $\pi _{1,t}(U)$ and $\pi _{2,u}
(U)$ are simply connected in $\bf C$ for each $t$ and $u\in
\bf C$ for which there exists $z\in U$, $U$ is an open subset in $\bf H$,
then there exists $g\in C^{\omega }(U,{\bf H})$ such that
$g'(z).I=f(z)$ for each $z\in U$.}
\par {\bf 3.11.3. Theorem.} {\it Let $U$ and $V$ be quaternion
holomoprhically simply connected regions in $\bf H$ with
$U\cap V\ne \emptyset $ connected.
Then $U\cup V$ is quaternion holomorphically simply connected.}
\par {\bf 3.12. Corollary.} {\it Let $U$ be an open subset in $\bf H^n$,
then the family of all quaternion holomorphic functions
$f: U\to \bf H$ has a structure of an $\bf H$-algebra.}
\par {\bf Proof.} If $f_1(z)=\alpha g(z)\beta +\gamma h(z)\delta $
or $f_2(z)=g(z)h(z)$ for each $z\in U$, where $\alpha $, $\beta $,
$\gamma $ and $\delta \in \bf H$ are constants, $g$ and $h$ are
quaternion holomorphic functions on $U$, then $F_1$ and $F_2$ are Frech\'et
differentiable on $U$ by $(v,w,x,y)$ (see \S 2.1 and \S 2.2)
and $D_{\tilde z}f_1(z)=\alpha (D_{\tilde z}g)\beta +\gamma
(D_{\tilde z}h)\delta =0$ and $D_{\tilde z}f_2(z)=
(D_{\tilde z}g)h+g(D_{\tilde z}h)=0$, hence $f_1$ and $f_2$ are also
quaternion holomorphic on $U$.
\par {\bf 3.13. Proposition.} {\it For each complex holomorphic function
$f$ in a neighbourhood $B(t_0,r,{\bf C})$
of a point $t_0\in \bf C$ there exists
a quaternion $z$-analytic function $g$ on a neighbourhood $B(a,r,{\bf H})$
of $a\in \bf H$ such that $a_{1,1}=t_0$
(or $a_{1,2}=t_0$) and $g_{1,1}(t,u_0)=f(t)$ (or $g_{1,2}(u_0,t)=
f(t)$ respectively)
on $B(t_0,r,{\bf C})$, where $B(x,r,X):=\{ y\in X: \rho _X(x,y)\le r \}$
is the ball in a space $X$ with a metric $\rho $, $r>0$, $u_0=a_{1,2}$
(or $u_0=a_{1,1}$ correspondingly).}
\par {\bf Proof.} Write Conditions $(2.2)$ for a right
superlinearly superdifferentiable function in the complex form.
This yields:
$$(3.12)\quad \partial f_{1,1}/\partial {\bar t}=0,
\quad \partial f_{1,2}/\partial t=0,
\quad \partial f_{1,1}/\partial u=0,
\quad \partial f_{1,2}/\partial {\bar u}=0.$$
There are also skew conditions:
$$(3.13)\quad \partial (g_{1,1}+ih_{1,2})/\partial (w+ix)=0,
\quad \partial (g_{1,2}+ih_{1,1})/\partial (w+ix)=0.$$
Other conditions derive from these.
For example, for the pair of variables $(v,y)$
using the matrix $L$ we get
$$(3.14)\quad \partial g_{1,1}/\partial v=\partial h_{1,2}/\partial y,
\quad \partial h_{1,2}/\partial v=- \partial g_{1,1}/\partial y,$$
$$(3.15)\quad \partial g_{1,2}/\partial v=\partial h_{1,1}/\partial y,
\quad \partial h_{1,1}/\partial v=- \partial g_{1,2}/\partial y,$$
which in complex form is the following:
$$(3.16)\quad \partial (g_{1,1}+ih_{1,2})/\partial
(v-iy)=0, \quad \partial (g_{1,2}+ih_{1,1})/\partial (v-iy)=0.$$
Equations $(3.13)$ and $(3.16)$ are equivalent to:
$$(3.13')\quad \partial (f_{1,1}+f_{1,2})/\partial (w+ix)=0,
\quad \partial ({\bar f}_{1,1}-{\bar f}_{1,2})/\partial (w+ix)=0,$$
$$(3.16')\quad \partial (f_{1,1}+f_{1,2})/\partial (v-iy)=0,
\quad \partial ({\bar f}_{1,1}-{\bar f}_{1,2})/\partial (v-iy)=0,$$
that is, there are two functions $p$ and $q$ holomorphic
in complex variables $w-ix$ and $v+iy$ such that
$f_{1,1}(z)=p(w-ix,v+iy)+q(w-ix,v+iy)$ and
$f_{1,2}(z)=p(w-ix,v+iy)-q(w-ix,v+iy)$.
\par Consider first an extension in the class
of quaternion holomorphic functions with a right superdifferential
not necessarily right superlinear on the superalgebra $\bf H^n$.
Since $f$ is holomorphic in $B(t_0,r,{\bf C})$, it has
a decomposition $f(t)=\sum_{n=0}^{\infty }f_n(t-t_0)^n$,
where $f_n\in \bf C$. Consider its extension in
$B({{t_0\quad 0}\choose {0\quad {\bar t}_0}},r,{\bf H})$
such that $f(z)=\sum_{n=0}^{\infty }
{{f_n\quad 0}\choose {0\quad {\bar f}_n}}
(z-{{t_0\quad 0}\choose {0\quad {\bar t}_0}})^n$.
Evidently this series converges for each $z\in
B({{t_0\quad 0}\choose {0\quad {\bar t}_0}},r,{\bf H})$
and this extension of $f$ is quaternion holomorphic, since
${{f_n\quad 0}\choose {0\quad {\bar f}_n}}\in \bf H$ for each $n$
and ${{t_0\quad 0}\choose {0\quad {\bar t}_0}}\in \bf H$, that is,
$\partial f/\partial {\tilde z}=0$.
If $z={{t\quad u}\choose {-{\bar u}\quad {\bar t}}}$
and $u=0$, then $f(z)=\sum_{n=0}^{\infty }
{{f(t)\quad 0}\choose {0\quad {\bar f}(t)}}$.
Another type of a solution is:
$f(z)=\sum_{n=0}^{\infty }{ {f_n\quad 0}\choose {0\quad {\bar f}_n} }$
$((z-JzJ)/2-{ {t_0\quad 0}\choose {0\quad {\bar t}_0 }} )^n$, since
${{t\quad 0}\choose {0\quad {\bar t} }}=(z-JzJ)/2$ for each
$z={{t\quad u}\choose {-{\bar u}\quad {\bar t} }}$.
\par Consider now more narrow class of quaternion holomorphic functions
with a right superlinear superdifferential on the superalgebra $\bf H^n$.
There is another way to construct $f$ on
$B({{t_0\quad 0}\choose {0\quad {\bar t}_0}},r,{\bf H})$, because
due to Theorems 2.15 and 3.10 a quaternion holomorphic function in
interior of this ball is quaternion $z$-analytic in it.
The construction of $f_{1,1}(t,u)$ satisfying the conditions above
and such that $f_{1,1}(t,0)=f(t)$ in $B(t_0,r,{\bf C})$ 
then comes down to
finding $f_{1,2}(t,u)$ with partial differential skew conditions
arising from those for $f_{1,1}$ and specific conditions on
$f_{1,2}$ such that $f_{1,1}$ is holomoprhic in $t$ and antiholomorphic in
$u$, $f_{1,2}$ is holomorphic in $u$ and antiholomorphic in
$t$ (where antiholomorphic means holomorphic in the
complex conjugate variable $\bar u$ or $\bar t$ respectively).
\par The second type of extension can be obtained from the first
by applying right multiplication by $\tilde K$ on the right, that is,
$f(z)\mapsto f(z){\tilde K}={{f_{1,2}\quad -f_{1,1}}\choose
{{\bar f}_{1,1}\quad {\bar f}_{1,2}}}.$
\par {\bf 3.14. Proposition.} {\it If $f$ is a quaternion
holomorphic function on an open subset $U$ in $\bf H$,
where $f'(z)\ne 0$ and $f'(z)$ is right superlinear,
then it is a conformal mapping in each point $z\in U$,
that is preserving angles between differentiable curves.}
\par {\bf Proof.} Let $z\in U$, then $f$ is differentiable
at $z$ and there exists $\lambda =f'(z)\in \bf H$.
Each quaternion $h={{h_t\quad h_u}\choose {-{\bar h}_u\quad
{\bar h}_t}} \in \bf H$ can be considered as vector
$(h_t,h_u)$ in $\bf C^2$.
Consider a scalar product in $\bf C^2$: $\quad $
$(h,k):=h_t{\bar k}_t+h_u{\bar k}_u$.
On the other hand, if $\lambda
\ne 0$, then $\lambda =|\lambda | \zeta $, where
$|\zeta |=1$. Rows and columns of the $2\times 2$ complex matrix
$\zeta $ are orthonormal, hence it is unitary and
$(\zeta h,\zeta k)=(h,k)$ for each $h$ and $k$ in $\bf C^2$
or for the corresponding quaternions in $\bf H$.
Therefore, for each vectors $h\ne 0$ and $k\ne 0$ in $\bf H$
$$(3.17)\quad (\lambda h,\lambda k)/(|\lambda h| |\lambda k|)=
(h,k)/(|h| |k|).$$
If $\psi $ and $\phi : (-1,1)\to U$ are two differentiable curves
crossing in a point $z\in U$, then we have two vectors $\psi '(0)=:h$
and $\phi '(0)=k$, where $\psi (0)=\phi (0)=z$.
Then $f(\psi (s))'=f'(z)|_{z=\psi (s)}.\psi '(s)$.
From Formula $(3.17)$ it follows, that $f$ preserves
the angle $\alpha $ between curves $\psi $ and $\phi $,
where $\cos (\alpha )=Re (\psi '(0),\phi '(0))/
(|\psi '(0)| |\phi '(0)|)$ for $\psi '(0)\ne 0$ and $\phi '(0)\ne 0$.
\par {\bf 3.15. Theorem.} {\it Let $f$ be a quaternion
holomorphic function on an open subset $U$ in $\bf H$ such that
$\sup_{z\in U, h\in B(0,1,{\bf H})} |[f(z)(\zeta -z)^{-2}]^{\hat .}.h|
\le C/|\zeta -z|^2$ for each $\zeta \in {\bf H}\setminus cl (U)$.
Then $|f'(z)|\le C/d(z)$ for each $z\in U$,
where $d(z):=\inf_{\zeta \in {\bf H}\setminus U}|\zeta -z|$;
$|f(\xi )-f(z)|/|\xi -z|\le 2C/r$ for each $\xi $ and $z\in
B(a,r/2,{\bf H})\subset Int (B(a,r,{\bf H}))\subset U$, where
$r>0$. In particular, if $f$ is a quaternion holomorphic 
function with bounded $[f(z)(\zeta -z)^{-2}]^{\hat .}.h|\zeta -z|^2$ on
${\bf H}^2\times B(0,1,{\bf H})$ with $|\zeta |\ge 2|z|$,
that is, $\sup_{\zeta ,z\in {\bf H}, |\zeta |\ge 2|z|,
h\in B(0,1,{\bf H})} |[f(z)(\zeta -z)^{-2}]^{\hat .}.h| |\zeta -z|^2
<\infty $, then $f$ is constant.}
\par {\bf Proof.} In view of Theorem $3.9$ there exists
a rectifiable path $\gamma $ in $U$ such that
$$(3.18)\quad (\partial ^kf(z)/\partial z^k)=
k!(2\pi )^{-1}(\int_{\gamma +z_0}f(\zeta )(\zeta -z)^{-k-1}d\zeta )
(wJ+xK+yL)^{-1},$$
where $\gamma (s)=r'\exp (2\pi s(wJ+xK+yL))$ with $s\in [0,1]$, $0<r'$.
Therefore, $|f'(z)|\le C/d(z)$. Since $\int_{\zeta }^zdf(z)=
f(z)-f(\zeta )$, then $|f(\xi )-f(z)|/|\xi -z|\le
\sup_{z\in B(a,r/2,{\bf H})}[C/d(z)]\le 2C/r$, where $r'<r/2$,
$\xi $ and $z\in B(a,r/2,{\bf H})\subset Int (B(a,r,{\bf H}))\subset U$.
Taking $r$ tending to infinity,
if $f$ is quaternion holomorphic with bounded
$[f(z)(\zeta -z)^{-2}]^{\hat .}.h |\zeta -z|^2$
on ${\bf H}^2\times B(0,1,{\bf H})$ for $|\zeta |\ge 2|z|$, then
$f'(z)=0$ for each $z\in \bf H$, since $f$ is locally $z$-analytic
and $\sup_{\zeta ,z\in U, |\zeta |\ge 2|z|,
h\in B(0,1,{\bf H})} |[f(z)(\zeta -z)^{-2}]^{\hat .}.h|
|\zeta -z|^2 <\infty $
is bounded, hence $f$ is constant on $\bf H$.
\par {\bf 3.16. Remark.} Theorems 3.9, 3.10 and 3.15 are
the quaternion analogs of the Cauchy, Morera and Liouville theorems
correspondingly. Evidently, Theorem 3.15 is also true for
right superlinear ${\hat f}(z)$ on $\bf H$ for each $z\in U$ and
with bounded ${\hat f}(z).h$ on $U\times B(0,1,{\bf H})$ instead of
$[f(z)(\zeta -z)^{-2}]^{\hat .}.h |\zeta -z|^2$. In particular, if
$f$ is quaternion holomorphic on $\bf H$ and ${\hat f}(z)$
is right superlinear on $\bf H$ for each $z\in \bf H$ and
${\hat f}(z).h$ is bounded on $U\times B(0,1,{\bf H})$,
then $f$ is constant.
\par {\bf 3.17. Theorem.} {\it Let $P(z)$ be a polynomial
on $\bf H$ such that $P(z)=z^{n+1}+\sum_{\eta (k)=0}^n(A_k,z^k)$,
where $A_k=(a_{1,k},...,a_{s,k}),$ $a_{j,l}\in \bf H$,
$k=(k_1,...,k_s)$, $0\le k_j\in \bf Z$, $\eta (k)=k_1+...+k_s$,
$0\le s=s(k)\in \bf Z$, $s(k)\le \eta (k)+1$,
$(A_k,z^k):=a_{1,k}z^{k_1}...a_{s,k}z^{k_s}$, $z^0:=1$.
Then $P(z)$ has a root in $\bf H$.}
\par {\bf Proof.} Consider a polynomial $Q(z):=z^{n+1}+\sum_{\eta (k)=0}^n
(z^k,{\tilde A}_k)$, where
$(z^k,{\tilde A_k}):=z^{k_1}{\tilde a}_{1,k}...z^{k_s}{\tilde a}_{s,k}.$
Then $PQ$ is a quaternion holomorphic function on $\bf H$.
Suppose that $P(z)\ne 0$ for each $z\in H$.
Consider a rectifiable path $\gamma _R$ in $\bf H$ such that
$\gamma _R([0,1])\cap {\bf H}=[-R,R]$ and outside $[-R,R]$:
$\quad \gamma _R(s)=R\exp (2\pi sM)$, where $M$ is a unit vector
in $\bf H_i$. Since $\lim_{|z|\to \infty }P(z)z^{-n-1}=1$, then
due to Theorem 2.11 $\lim_{R\to \infty }\int_{\gamma _R}(PQ)^{-1}(z)dz=$
$\int_{-R}^R(PQ)^{-1}(v)dv$ $=\int_{-R}^R|P(v)|^{-2}dv\ge 0.$
The latter integral is equal to zero if and only if $|P(v)|^{-2}=0$
for each $v\in \bf R$. This contradicts our supposition, hence
there exists a root $z_0\in \bf H$, that is, $P(z_0)=0$.
\par {\bf 3.17.1. Note.} Consider, for example,
the polynomial $P(z)=(z-a)^2+J(z-a)K$ on $\bf H$, then there does not exist
$\lim_{z\to a, z\ne a} f(z)(z-a)^{-1}$ and there also does not
exist $\lim_{z\to a, z\ne a} (z-a)^{-1}f(z)$, though $f(a)=0$.
This makes an obstacle for a quaternion analog of the Gauss
theorem about zeros of a derivative of a complex polynomial.
Even in a particular case, when a polynomial has a decomposition
$f(z)=(z-a_1)...(z-a_m)$, where $a_1$,...,$a_m\in \bf H$,
then $f(z)^{-1}f'(z).h=\sum_j (z-a_m)^{-1}...(z-a_{j+1})^{-1}(z-a_j)^{-1}
h(z-a_{j+1})...(z-a_m)$, consequently,
$(f(z)^{-1}f'(z).I)^{\tilde .}$ $=\sum_j\lambda _j(z-a_j)
\lambda _j^{-1}/|z-a_j|^2$, where $\lambda _j=
[(z-a_{j+1})...(z-a_m)]^{\tilde .}$. Hence $z\sum_jc_j|z-a_j|^{-2}=
\sum_j\lambda _ja_j\lambda _j^{-1}$, where $\lambda _jz=:
zc_j\lambda _j$, $c_j\in \bf H$, $|c_j|=1$, $z$ is a root of $f'(z).I$.
In the case of pairwise commuting $a_1$,...,$a_m$ this formula simplifies,
but in general $a_1$,...,$a_m$ need not be commuting.
\par {\bf 3.17.2. Remark.} The noncommutative geometry
in terms of a scheme theory for associative algebras
depends heavily on sheaf theory own a Zariski topology, even
a noncommutative version thereof \cite{oystaey}.
In this theory, instead of
starting from a noncommutative algebra and dealing with its
geometry as being "virtual" we now can consider concretely
defined geometrical objects, but defined over a noncommutative
field $\bf H$. The noncommutative $\bf H$-algebras that appear
are the rings of locally analytic functions in an open set
$U$ for the real topology, i.e. $C^{\omega }_{(z,\tilde z)}(U,{\bf H})$.
This obviously leads to the possibility to define presheaves
and sheaves on the objects embedded into $\bf H^n$ and endowed with
the induced real topology; it also applies to $\bf H$-analytic objects
like the four dimensional manifold $W$ constructed before
Theorem 3.8, and indeed to any quaternion version of a general
manifold, that is a "manifold" with a local $\bf H^n$-structure
generalizing in the obvious way the local $\bf R^m$-structure.
In later work we aim to study the quaternion version
of sheaf cohomology and Cartan Theorems A and B, as well
as noncommutative Stein manifolds, i.e. the quaternion
version of holomorphy domains.
\par {\bf 3.18. Theorem.} {\it Let $f$ be a quaternion holomorphic
function on an open subset $U$ in $\bf H$. Suppose that $\epsilon >0$
and $\sf K$ is a compact subset of $U$. Then there exists a function
$g(z)=P_{\infty }(z)+\sum_{k=1}^{\nu }P_k[(z-a_k)^{-1}]$,
$z\in {\bf H}\setminus \{ a_1,...,a_{\nu } \} $, $\nu \in \bf N$,
where $P_{\infty }$ and $P_j$ are polynomials, $a_j\in Fr (U)$,
$Fr (U)$ denotes a topological boundary of $U$ in $\bf H$,
such that $|f(z)-g(z)|<\epsilon $ for each $z\in \sf K$.}
\par {\bf Proof} is analogous to the proof of Runge's theorem
(see \cite{heins}) due to Theorem 3.9 and considering
four dimensional cubes $S_{j,k}=\mbox{ }_1S_j\times \mbox{ }_2S_k$
with ribs of length
$n^{-1}$ in $\bf H$ instead of two dimensional cubes
in $\bf C$ and putting $S:=\cup_{j,k}S_{j,k}$ such that
${\sf K}\subset Int (S)$, where $n\in \bf N$ tends to infinity,
$\mbox{ }_1S_j$ and $\mbox{ }_2S_k$
are two dimensional cubes in $\bf C_1$ and $\bf C_2$
which are two copies of $\bf C$ embedded orthogonally in $\bf H$
as ${\bf R}I\oplus {\bf R}J$ and ${\bf R}K\oplus {\bf R}L$
correspondingly.
Since $f$ is quaternion holomorphic and $\bf  K$ is compact,
we may apply Formula $(3.9)$ to each $B_{j,k}=\gamma $
such that $\gamma ={{\gamma _{1,1}\quad \gamma _{1,2}}\choose
{-{\bar \gamma }_{1,2}\quad {\bar \gamma }_{1,1}}}$,
$\gamma _{1,1}=\partial \mbox{ }_1S_j,$ $\gamma _{1,2}=\partial
\mbox{ }_2S_k$, it can be seen, that $f$ can be approximated uniformly on $\sf K$
by a sum of the form $\sum_{k=1}^{\mu }(a_{1,k}(\zeta _k-z)^{-1}
a_{2,k})$, where $a_{j,k}\in \bf H$,
$\zeta _k\in Fr (S)$. 
For a given $n\in \bf N$ if $b\in Fr (S)$, then there exists
$a\in Fr (U)\cup \partial B(0,r,{\bf H})$ such that
$|b-a|\le n^{-1}$. If $z\in \sf K$ and $|z-a|\ge n^{-1}$, then
the series $(z-b)^{-1}=(\sum_{k=0}^{\infty }[(z-a)^{-1}(b-a)]^k)
(z-a)^{-1}$ converges uniformly on $\sf K$ and it is clear that
$f$ can be approximated uniformly on $\sf K$ by a function of
the indicated form (see also \S 3.17).
\par {\bf 3.19. Note and Definitions.} Consider a one-point
(Alexandroff) compactification $\hat {\bf H}$ of the locally compact
topological space $\bf H$. It is homeomorphic to a unit
four dimensional sphere $S^4$ in the Euclidean space
$\bf R^5$. If $\zeta $ is a point in $S^4$ different
from $(1,0,0,0,0)$, then the straight line containing
$(1,0,0,0,0)$ and $\zeta $ crosses $\pi _S$ in a finite point
$z$, where $\pi _S$ is the four dimensional plane orthogonal to the vector
$(1,0,0,0,0)$ and tangent to $S^4$ at the south pole
$(-1,0,0,0,0)$. This defines the bijective continuous mapping from
$S^4\setminus \{ (1,0,0,0,0) \} $ onto
$\pi _S$ such that $(1,0,0,0,0)$ corresponds to
the point of infinity. Therefore each function on a subset $U$
of $\bf H$ as a topological space
can be considered on the homeomorphic subset $V$ in
$S^4$. If $U$ is a locally compact subset of
$\bf H$ and $\lim_{z\in U, |z|\to \infty }
f(z)$ exists, then $f$ has an extension on $\hat U$.
\par Let $z_0\in \hat {\bf H}$ be a marked point.
If a function $f$ is defined and quaternion holomorphic 
on $V\setminus \{ z_0 \} $,
where $V$ is a neighbourhood of $z_0$, then $z_0$ is called
a point of an isolated singularity of $f$.
\par Suppose that $f$ is a quaternion holomorphic function
in $B(a,0,r,{\bf H})\setminus \{ a \} $ for some $r>0$.
Then we say that $f$ has an isolated singularity at $a$.
Let $B(\infty ,r,{\bf H}):=\{ z\in {\hat {\bf H}}$ such that
$r^{-1}<|z|\le \infty \} $. The we say that $f$ has an isolated
singularity at $\infty $ if it is quaternion holomorphic
in some $B(\infty ,r,{\bf H})$.
\par Let $f: U\to \bf H$ be a function, where $U$ is a neighbourhood
of $z\in \hat {\bf H}$. Then $f$ is said to be meromorphic at
$z$ if $f$ has an isolated singularity at $z$.
If $U$ is an open subset in $\hat {\bf H}$, then $f$ is called
meromorphic in $U$ if $f$ is meromorphic at each point $z\in U$.
If $U$ is a domain of $f$ and $f$ is meromorphic in $U$, then
$f$ is called meromorphic on $U$. Denote by ${\bf M}(U)$ the set
of all meromorphic functions on $U$.
Let $f$ be meromorphic on
a region $U$ in $\hat {\bf H}$. A point $c\in \bigcap_{V\subset
U, V\mbox{ is compact }}cl( f(U\setminus V))$ is called a cluster
value of $f$.
\par {\bf 3.20. Proposition.} {\it Let $f$ be a
quaternion holomorphic function with a right
$\bf H$-superlinear superdifferential on
an open connected subset $U\subset \hat {\bf H}$ and
suppose that there exists
a sequence of points $z_n\in U$ having a cluster point $z\in U$
such that $f(z_n)=0$ for each $n\in \bf N$,
then $f=0$ everywhere on $U$.}
\par {\bf Proof} follows from the local $z$-analyticity of
$f$ and the fact $f^{(k)}(z)=0$ for each $0\le k\in \bf Z$
(see Theorems 2.11 and 3.10), when $f'(z)$ is right
$\bf H$-superlinear on $U$, since  \\
$f^{(k)}(z)=\lim_{n+m\to \infty }
(f^{(k-1)}(\zeta _n)- f^{(k-1)}(\zeta _m))
(\zeta _n-\zeta _m)^{-1},$
where $\zeta _n$ is a subsequence of $\{ z_n: n \} $ of pairwise
distinct points converging to $z$. Therefore, $f$ is equal to
zero on a neighbourhood of $z$. The maximal subset of $U$
on which $f$ is equal to zero is open in $U$. On the other hand it is closed,
since $f$ is continuous, hence $f$ is equal to zero on $U$, since
$U$ is connected.
\par {\bf 3.21. Note.} Without the condition of right $\bf H$-superlinearity
of $f'(z)$ on $U$ Proposition 3.20 is not true in general,
since $f_1(z):=azb$ and $f_2(z):=abz$ coincide on ${\bf R}I$,
but not on any neighbourhood of zero, when $a$ and $b$
are noncommuting fixed quaternions, $z\in \bf H$.
\par Consider a function $f(z)=z^{-1}az$
on ${\bf H}\setminus \{ 0 \} $,
where $0\ne z\in \bf H$ and $a\in \bf H$.
If $a\ne vI$, then there exists $0\ne h\in \bf H$
such that $h^{-1}ah=:b\ne a$. For $z=sh$ there exists
$\lim_{z=sh, s\ne 0, s\to 0}f(z)=b$, for $z=sI$
$\lim_{z=sI, s\ne 0, s\to 0}f(z)=a$.
Therefore, if $a\ne vI$ for some $v\in \bf R$, then there does not exist
a limit of $f(z)$ for $z$ tending to zero.
This makes clear, that for the quaternion field
it is important to consider an analog of a Laurent series
of a function quaternion holomorphic on $U\setminus \{ 0 \} $
not only in terms $az^k$, but also in $a_1z^{k_1}...a_nz^{k_n}$,
where $k_j$ are integers, $z^0:=1$.
\par {\bf 3.22. Theorem.} {\it Let $\bf A$ denote the family of
all functions $f$ such that $f$ is quaternion holomorphic
on $U:=Int (B(a,r,R,{\bf H}))$,
where $a$ is a marked point in $\bf H$, $0\le r<R<\infty $
are fixed. Let $\bf S$ denote a subset
of ${\bf Z}^{\bf N}$ such that for each $k\in S$ there exists
$m(k):=\max \{ j: $ $k_j\ne 0, k_i=0 $ $\mbox{for each }
i>j \} \in \bf N$ and let ${\bf B}$ be a family of finite sequences
$b=(b_1,...,b_n)$ such that $b_j\in \bf H$ for each $j=1,...,n$,
$n\in \bf N$.
Then there exists a bijective correspondence between
$\bf A$ and $q\in {\bf B}^S$ such that
$$(3.19) \quad \lim_{m+\eta
\to \infty }\sup_{z\in B(a,r_1,R_1,{\bf H})}
\sum_{k, m(k)=m, \eta (k)=\eta } |(b_k,z^k)|=0$$
for each $r_1$ and $R_1$ such that $r<r_1<R_1<R$,
where $\eta (k):=k_1+...+k_{m(k)}$, $q(k)=:b_k
=(b_{k,1},...,b_{k,m(k)})$,
$(b_k,z^k)=b_{k,1}z^{k_1}...b_{k,m(k)}z^{k_{m(k)}}$
for each $k\in \bf S$, that is,
$f\in \bf A$ can be presented by a convergent series
$$(3.20)\quad f(z)=\sum_{b\in q}(b_k,z^k).$$}
\par {\bf Proof.} If Condition $(3.19)$ is satisfied, then
the series $(3.20)$ converges on $B(a,r',R',{\bf H})$ for each
$r'$ and $R'$ such that $r<r'<R'<R$, since $r_1$ and $R_1$ are
arbitrary such that $r<r_1<R_1<R$ and $\sum_{n=0}^{\infty }p^n$
converges for each $|p|<1$. In particular taking
$r_1<r'<R'<R_1$ for $p=R'/R_1$ or
$p=r_1/r'$. Therefore, from $(3.19)$ and
$(3.20)$ it follows, that $f$ presented by the series $(3.20)$
is quaternion holomorphic on $U$.
\par Vice versa let $f$ be in $\bf A$. In view of Theorems 2.11 and 3.9
there are two rectifiable closed paths $\gamma _1$ and $\gamma _2$
such that $\gamma _2(s)=a+r' \exp (2\pi s M_2)$ and
$\gamma _1(s)=a+R' \exp (2\pi s M_1)$, where $s\in [0,1]$,
$M_1$ and $M_2\in \bf H$ with $|M_1|=1$ and $|M_2|=1$,
where $r<r'<R'<R$, because as in \S 3.9 $U$ can be presented
as a finite union of regions $U_j$ each of which satisfies
the conditions of Theorem 2.11. Using a finite number
of rectifiable paths $w_j$ (joining $\gamma _1$ and $\gamma _2$
within $U_j$) going twice in one and the
opposite directions leads to the conclusion
that for each $z\in Int (B(a,r',R',{\bf H}))$
the function $f(z)$ is presented by the integral formula:
$$(3.21)\quad f(z)=(2\pi )^{-1}\{
(\int_{\gamma _1}f(\zeta )(\zeta -z)^{-1} d\zeta )M_1^{-1}-
(\int_{\gamma _2}f(\zeta )(\zeta -z)^{-1} d\zeta )M_2^{-1} \} .$$
On $\gamma _1$ we have the inequality:
$|(\zeta -a)^{-1}(z-a)|<1$, on $\gamma _2$ another
inequality holds: $|(\zeta -a)(z-a)^{-1}|<1$.
Therefore, for $\gamma _1$ the series
$$(\zeta -z)^{-1}=(\sum_{k=0}^{\infty }((\zeta -a)^{-1}(z-a))^k)
(\zeta -a)^{-1}$$
converges uniformly by $\zeta \in B(a,R_2+\epsilon,R_1,{\bf H})$ and
$z\in B(a,r_2,R_2,{\bf H})$, while for $\gamma _2$ the series
$$(\zeta -z)^{-1}=-(z-a)^{-1}(\sum_{k=0}^{\infty }((\zeta -a)(z-a)^{-1})^k)$$
converges uniformly by $\zeta \in B(a,r_1,r_2-\epsilon,{\bf H})$
and $z\in B(a,r_2,R_2,{\bf H})$ for each $r'<r_2<R_2<R'$ and each
$0<\epsilon <\min (r_2-r_1,R_1-R_2)$. Consequently,
$$(3.22)\quad f(z)=\sum_{k=0}^{\infty }(\phi _k(z)+\psi _k(z)),
\mbox{ where}$$
$$\phi _k(z):=(2\pi )^{-1} \{ \int_{\gamma _1}
f(\zeta )((\zeta -a)^{-1}(z-a))^k(\zeta -a)^{-1} d\zeta )M_1^{-1} \} ,$$
$$\psi _k(z):=(2\pi )^{-1}\{ \int_{\gamma _2}
f(\zeta )(z-a)^{-1}((\zeta -a)(z-a)^{-1})^k d\zeta )M_2^{-1} \} ,$$
and where $\phi _k(z)$ and $\psi _k(z)$ are quaternion holomorphic functions,
hence $f$ has decomposition $(3.20)$ in $U$, since due to \S 2.15
and \S 3.9 there exists $\delta >0$ such that integrals for $\phi _k$
and $\psi _k$ by $\gamma _1$ and $\gamma _2$ are the same for each
$r'\in (r_1,r_1+\delta )$, $R'\in (R_1-\delta ,R_1)$.
Using the definition of the quaternion line integral we get $(3.20)$
converging on $U$. Varying $z\in U$ by $|z|$ and $Arg (z)$ we get
that $(3.20)$ converges absolutely on $U$, consequently,
$(3.19)$ is satisfied.
\par {\bf 3.23. Notes and Definitions.} Let $\gamma $ be a closed
curve in $\bf H$. There are natural projections from $\bf H$
on complex planes: $\pi _1(z)=(v,w)$, $\pi _2(z)=(v,x)$,
$\pi _3(z)=(v,y)$, where $z=vI+wJ+xK+yL$ with real $v$, $w$, $x$
and $y$. Therefore, $\pi _l(\gamma )=:\gamma _l$ are curves
in complex planes $\bf C_1$ isomorphic to ${\bf R}I\oplus {\bf R}J$,
$\bf C_2$ isomorphic to ${\bf R}I\oplus {\bf R}K$ and
$\bf C_3$ isomorphic to ${\bf R}I\oplus {\bf R}L$, where
$l=1,2,3$ respectively. A curve $\gamma $ in $\bf H$ is closed
(a loop, in another words) if and only if $\gamma _l$
are closed for each $l=1,2,3$, that is, $\gamma (0)=\gamma (1)$
and $\gamma _l(0)=\gamma _l(1)$ correspondingly.
For each point $a\in \bf H$ we have its projections $a_l:=\pi _l(a)$.
In each complex plane there is the standard complex notion of
a topological index $In (a_l,\gamma _l)$ of a curve $\gamma _l$ at $a_l$.
Therefore, there exists a vector $In (a,\gamma ):=\{ In (a_1,\gamma _1),
In (a_2,\gamma _2), In (a_3,\gamma _3) \} $ which we call
the topological index of $\gamma $ at a point $a\in \bf H$.
This topological index is invariant relative to homotopies satisfying
conditions of Theorem 3.9.
Consider now a standard closed curve $\gamma (s)=a+r\exp (2\pi snM)$,
where $M\in \bf H_i$
with $|M|=1$, $n\in \bf Z$, $r>0$, $s\in [0,1]$.
Then ${\hat I}n (a,\gamma ):=(2\pi )^{-1}(\int_{\gamma }d Ln (z-a))=nM$
is called the quaternion index of $\gamma $ at a point $a$.
It is also invariant relative to homotopies satisfying the
conditions of Theorem
3.9. Moreover, ${\hat I}n (h_1ah_2,h_1\gamma h_2)={\hat I}n (a,\gamma )$
for each $h_1$ and $h_2\in {\bf H}\setminus \{ 0 \} $
such that $h_1Mh_2=M$.
For $M=wJ+xK+yL$ there is the equality ${\hat I}n (a,\gamma )=
In (a_1,\gamma _1)wJ+ In (a_2,\gamma _2)xK+In (a_3,\gamma _3)yL$
(adopting the corresponding convention for signs of indexes
in each $\bf C_j$ and the convention of positive directions of going along
curves).
In view of the properties of $Ln$ for each curve $\psi $ in $\bf H$
there exists $\int_{\gamma }d Ln (z-a)=2\pi qM$ for some
$q\in \bf R$ and $M\in \bf H_i$ with $|M|=1$.
For a closed  curve $\psi $ up to a composition of homotopies
each of which is charaterized by homotopies in $\bf C_l$
for $l=1,2,3$ there exists a standard $\gamma $ with a generator $M$
for which ${\hat I}n (a,\gamma )=qM$, where $q\in \bf Z$.
Therefore, we can take as a definiton
${\hat I}n (a,\psi )={\hat I}n (a,\gamma )$.
Define also the residue of a meromorphic function with an isolated
singularity at a point $a\in \bf H$ as
$res (a,f):=(\int_{\gamma }f(z)dz)(2\pi M)^{-1}$,
where $\gamma (s)=a+r\exp (2\pi sM)\subset V,$ $r>0$,
$|M|=1$, $M\in \bf H_i$, $s\in [0,1]$,
$f$ is quaternion holomorphic on $V\setminus \{ a \} $.
\par If $f$ has an isolated singularity at $a\in {\hat {\bf H}}$,
then coefficients $b_k$ of its Laurent series (see \S 3.22)
are independent of $r>0$. The common series is called
the $a$-Laurent series. If $a=\infty $, then $g(z):=f(z^{-1})$
has a $0$-Laurent series $c_k$ such that $c_{-k}=b_k$.
Let $\beta := \sup_{b_k\ne 0}
\eta (k)$, where $\eta (k)=k_1+...+k_m$, $m=m(k)$
for $a=\infty $; $\beta =\inf_{b_k\ne 0} \eta (k)$ for $a\ne \infty $.
We say that $f$ has a removable singularity, pole,
essential singularity at $\infty $ according as
$\beta \le 0$, $0<\beta <\infty $, $\beta =+\infty $.
In the second case $\beta $ is called the order of the pole
at $\infty $. For a finite $a$ the corresponding cases
are: $\beta \ge 0$, $-\infty <\beta <0$, $\beta =-\infty $.
If $f$ has a pole at $a$, then $|\beta |$ is called
the order of the pole at $a$.
\par A value of a function $\partial _f(a):=\inf \{ \eta (k): b_k\ne 0 \} $
is called a divisor of $f$ at $a\ne \infty $,
$\partial _f(a):=\inf \{ - \eta (k): b_k\ne 0 \} $ for $a=\infty $.
Then $\partial _{f+g}(a)\ge \min \{ \partial _f(a), \partial _g(a) \} $
for each $a\in dom (f)\cap dom (g)$
and $\partial _{fg}(a)=\partial _f(a)+\partial _g(a).$
For a function $f$ meromorphic
on an open subset $U$ in $\bf \hat H$ the function $\partial _f(p)$
by the variable $p\in U$ is called the divisor of $f$.
\par {\bf 3.24. Theorem.} {\it Let $U$ be an open region in $\bf \hat H$
with $n$ distinct marked points $p_1,...,p_n$, and let $f$ be a quaternion
holomorphic function on $U\setminus \{ p_1,...,p_n \} =:U_0$
and $\psi $ be a rectifiable closed curve lying in
$U_0$ such that $U_0$ satisfies the conditions of Theorem 3.9
for each $z_0\in \{ p_1,...,p_n \} $. Then
$$\int_{\gamma }f(z)dz=2\pi \sum_{j=1}^n {\hat I}n (p_j,\gamma )
res (p_j,f).$$}
\par {\bf Proof.} For each $p_j$ consider the principal part $T_j$
of a Laurent series for $f$ in a neighbourhood of $p_j$,
that is, $T_j(z)=\sum_{k, \eta (k)<0} (b_k,(z-p_j)^k)$,
where $\eta (k)=k_1+...+k_n$ for $k=(k_1,...,k_n)$
(see Theorem 3.22).
Therefore, $h(z):=f(z)-\sum_jT_j(z)$ is a function having
a quaternion holomorphic extension on $U$.
In view of Theorem 3.9 for a quaternion holomorphic
function $g$ in a neighbourhood $V$ of a point $p$ and a rectifiable
closed curve $\zeta $ we have
$${\hat I}n (p,\zeta )g(p)=(2\pi )^{-1}(\int_{\zeta }g(z)(z-p)^{-1}dz)$$
(see \S 3.23).
We may consider small $\zeta _j$ around each $p_j$ with
${\hat I}n (p_j,\zeta _j)={\hat I}n (p_j,\gamma )$
for each $j=1,...,n$.
Then $\int_{\zeta _j}f(z)dz=\int_{\zeta _j}T_j(z)dz$ for each $j$.
Representing $U_0$ as a finite union of open regions $U_j$
and joining $\zeta _j$ with $\gamma $ by paths $\omega _j$ going
in one and the opposite direction as in Theorem 3.9 we get
$$\int_{\gamma }f(z)dz+\sum_j\int_{\zeta _j^-}f(z)dz=0,$$
consequently,
$$\int_{\gamma }f(z)dz=\sum_j\int_{\zeta _j}f(z)dz=
\sum_j 2\pi {\hat I}n (p_j,\gamma ) res (p_j,f),$$
where ${\hat I}n(p_j,\gamma )$ and $res (p_j,f)$ are
invariant relative to homotopies satisfying conditions of
Theorem 3.9.
\par {\bf 3.25. Corollary.} {\it Let $f$ and $T$
be the same as in \S 3.24, then $res (p_j,f)=res (p_j,T_j)=
res (p_j,\sum_{k,\eta (k)=-1}(b_k,(z-p_j)^k))$,
in particular, $res (p_j,b (z-p_j)^{-1})=b$.}
\par {\bf 3.26. Corollary.} {\it Let $U$ be an open region
in $\bf \hat H$ with $n$ distinct points $p_1,...,p_n$,
let also $f$ be a quaternion holomorphic function
on $U\setminus \{ p_1,...,p_n \} =:U_0$, $p_n=\infty $,
and $U_0$ satisfies conditions of Theorem 3.9 with at least
one $\psi $, $\gamma $ and each $z_0\in \{ p_1,...,p_n \} $.
Then $\sum_{p_j\in U} res (p_j,f)=0$.}
\par {\bf Proof.} If $\gamma $ is a closed curve encompassing
$p_1$,...,$p_{n-1}$, then $\gamma ^-(s):=\gamma (1-s)$, where $s\in [0,1]$,
encompasses $p_n=\infty $ with positive going by $\gamma ^-$ relative to
$p_n$. Since $\int_{\gamma }f(z)dz+\int_{\gamma ^-}f(z)dz=0$,
we get the satement of this corollary from Theorem 3.24.
\par {\bf 3.27. Definitions.} Let $f$ be a holomorphic function
on a neighbourhood $V$ of a point $z\in \bf H$.
Then the infimum: $n(z;f):=\inf \{ k: k\in {\bf N},
f^{(k)}(z)\ne 0 \} $ is called a multiplicity of $f$ at $z$.
Let $f$ be a holomorphic function on an open subset $U$
in $\bf \hat H$. Suppose $w\in \bf \hat H$,
then the valence $\nu _f(w)$ of $f$ at $w$ is by the definition
$\nu _f(w):=\infty $, when the set $\{ z: f(z)=w \} $
is infinite, and otherwise $\nu _f(w):=\sum_{z, f(z)=w}n(z;f)$.
\par {\bf 3.27.1. Theorem.} {\it Let $f$ be a meromorphic function
on a region $U\subset \bf \hat H$. If $b\in \bf \hat H$
and $\nu _f(b)<\infty $, then $b$ is not a cluster value of $f$
and the set $ \{ z: \nu _f(z)=\nu _f(b) \} $ is a neighbourhood
of $b$. If $U\ne \bf \hat H$ or $f$ is not constant, then
the converse statement holds. Nevertheless, it is false, when
$f=const $ on $\bf \hat H$.}
\par {\bf 3.27.2. Theorem.} {\it Let $U$ be a proper open subset
of $\bf \hat H$, let also $f$ and $g$ be two continuous functions
from ${\bar U}:=cl (U)$ into $\bf \hat H$ such that on a topological boundary
$Fr (U)$ of $U$ they satisfy the inequality
$|f(z)|<|g(z)|$ for each $z\in Fr (U)$. Suppose $f$ and $g$ are meromorphic
functions in $U$ and $h$ be a unique continuous map from
$\bar U$ into $\bf \hat H$ such that $h|_E=f|_E +g|_E$, where
$E:= \{ z: f(z)\ne \infty , g(z)\ne \infty \} $. Then
$\nu _{g|_U}(0)-\nu _{g|_U}(\infty )=\nu _{h|_U}(0)-\nu _{h|_U}(\infty )$.}
\par {\bf Proofs} of these two theorems are analogous to that of
Theorems VI.4.1, 4.2 \cite{heins}.
\par {\bf 3.28. Theorem.} {\it Let $U$ be an open subset
in $\bf H^n$, then there exists a represenation
of the $\bf R$-linear space
$C^{\omega }_{z,\tilde z}(U,{\bf H})$ of locally $(z,\tilde z)$-analytic
functions on $U$ such that it is isomorphic to the
$\bf R$-linear space $C^{\omega }_z(U,{\bf H})$
of quaternion holomorphic functions on $U$.}
\par {\bf Proof.} Evidently, the proof can be reduced to the case
$n=1$ by induction considering local $(z,\tilde z)$-series
decompositions by $(\mbox{ }^nz,\mbox{ }^n{\tilde z})$
with coefficients being convergent series of $(\mbox{ }^1z,
\mbox{ }^1{\tilde z},...,\mbox{ }^{n-1}z,\mbox{ }^{n-1}{\tilde z})$.
For each $z\in \bf H$ there are identities:
$JzJ=-vI-wJ+xK+yL$, $KzK=-vI+wJ-xK+yL$, $LzL=-vI+wJ+xK-yL$,
where $z=vI+wJ+xK+yL$ with $v$, $w$, $x$ and $y\in \bf R$.
Hence $P_{{\bf R}I}(z)=vI=(z-JzJ-KzK-LzL)/4$,
$P_{{\bf R}J}(z)=wJ=(z-JzJ+KzK+LzL)/4$,
$P_{{\bf R}K}(z)=xK=(z+JzJ-KzK+LzL)/4$,
$P_{{\bf R}L}(z)=yL=(z+JzJ+KzK-LzL)/4$ are projection operators
on ${\bf R}I$, ${\bf R}J$, ${\bf R}K$ and ${\bf R}L$ respectively,
where $I$, $J$, $K$ and $L$ are orthogonal vectors relative to the scalar
product in $\bf C^4$, ${\bf H}\ni z\mapsto (t,u,-{\bar u},{\bar t})
\in \bf C^4$.
Therefore, ${\tilde z}=vI-wJ-xK-yL=-(z+JzJ+KzK+LzL)/2$
and $d{\tilde z}=(dv)I-(dw)J-(dx)K-(dy)L=-(dz+J(dz)J+K(dz)K+L(dz)L)/2$.
Consequently, each polynomial in $(z,\tilde z)$ is also a polynomial
in $z$ only, moreover, each polynomial locally $(z,\tilde z)$
analytic function on $U$ is polynomial locally $z$-analytic on $U$.
Then if a series by $(z,{\tilde z})$ converges in a ball
$B(z_0,r,{\bf H^n})$, then its series in the $z$-representation
converges in a ball $B(z_0,r/2,{\bf H^n})$.
Then $\int_{\gamma }{\tilde z}dz=$ $-(\int_{\gamma }zdz+\int_{\gamma}
JzJdz$ $+\int_{\gamma }KzKdz+\int_{\gamma }LzLdz)/2$
and $\int_{\gamma }{\tilde z}dz=0$ for a closed rectifiable curve
$\gamma $ in $\bf H$ in such representation. This is not
contradictory, because from
$f_1|_{\gamma }=f_2|_{\gamma }$ it does not follow ${\hat f}_1|_{\gamma }
={\hat f}_2|_{\gamma }$, since ${\hat f}(z)$ is defined by values
of a function $f$ on an open neighbourhood of a point $z\in \bf H$,
where $f$, $f_1$ and $f_2\in C^0(U,{\bf H})$.
Therefore, $\int_{\gamma }d Ln z$ is quite different in general from
$\int_{\gamma }{\tilde z}dz$ (see \S 2.5 and \S 3.8). 
Considering basic polynomials
of any polynomial basis in $C^{\omega }_{z,\tilde z}(U,{\bf H})$
we get (due to infinite dimensionality of this space) a polynomial base of
$C^{\omega }_z(U,{\bf H})$. This establishes the 
$\bf R$-linear isomorphism between these two spaces.
Moreover, in such representation of the space
$C^{\omega }_{z,\tilde z}(U,{\bf H})$ we can put $D_{\tilde z}=0$,
yielding for differential forms $\partial _{\tilde z}=0$,
this leads to differential calculus and integration with
respect to $D_z$ and $dz$ only.
\par {\bf 3.29. Notes.} The latter paragraph also shows that for
$\mbox{ }_lC^{\omega }_{z,\tilde z}(U,{\bf H})$ and for
$\mbox{ }_rC^{\omega }_{z,\tilde z}(U,{\bf H})$ operators
$D_z$ and $D_{\tilde z}$ are different and neither $D_z$ nor
$D_{\tilde z}$ may be omitted from the differential calculus, since
automorphisms $z\mapsto azb$ of $\bf H$ with given quaternions
$a$ and $b$ such that $ab\ne 0$ do not leave 
$\mbox{ }_lC^{\omega }_{z,\tilde z}(U,{\bf H})$ and 
$\mbox{ }_rC^{\omega }_{z,\tilde z}(U,{\bf H})$ invariant.
\par Apart from the complex polynomial case in the quaternion case
a polynomial may have infinite family of roots, for example,
$P(z)=z^2+zJzJ+zKzK+zLzL-1$ has a $3$-dimensional over $\bf R$
manifold of roots $P(z)=0$, since $P(z)=-2|z|^2-1$.
\par Theorem 3.27.2 is the quaternion analog of the Rouch\'e theorem.
\par The function $f(z):=\cos (z{\tilde z}):=[\exp (Jz{\tilde z})+
\exp (-Jz{\tilde z})]/2$ is bounded on $\bf H$,
but neither the operator ${\hat f}$ is right superlinear, nor
the operator
$[f(z)(\zeta -z)]^{-2}]^{\hat .}.h|\zeta -z|^2$ is bounded on
${\bf H}^2\times B(0,1,{\bf H})$ for $|\zeta |\ge 2|z|$,
then $f(z)=\cos (z{\tilde z})$
can be written in the corresponding representation
as a quaternion locally $z$-analytic function, since $\sum_{n=1}^{\infty }
(2R)^n/n!$ converges for each $0\le R<\infty $.
This shows that in the last part
of Theorem 3.15 its conditions cannot be replaced by boundedness
of a quaternion holomorphic function $f$.
An interesting analog of the Liouville theorem for real
harmonic functions was investigated in \cite{mitid}.
Possibly the particular case of the quaternion analog
of the Liouville theorem for right superlinear $\hat f$
may be deduced from \cite{mitid} with the help
of Equations $(2.4)$ of Corollary 2.4 above.
\par There are other ways to define superdifferentiations
of algebras of quaternion functions:
\par $(1)$ factorize an algebra of quaternion locally
$(z,\tilde z)$-analytic functions $f: U\to \bf H$ by all relations
of the form $[\sum_j S_{j,1}zS_{j,2}-{\tilde z}]$, where
$S_{j,k}\in \bf H$ are fixed and $\sum_jS_{j,1}zS_{j,2}=\tilde z$
for each $z\in \bf  H$;
\par $(2)$ use as a starting point superlinearly superdifferentiable
functions $f: U\to \bf H$ and then prolong a superdifferentiation
on products of such functions with milder conditions on a
superdifferential, but they lead to the same result.
This approach can be generalized for general Clifford algebras
over $\bf R$, but some results will be weaker
or take another form, than in the case of $\bf H^n$.
\par {\bf 3.30. Theorem} (Argument principle). {\it
Let $f$ be a quaternion holomorphic function on an open region
$U$ satisfying conditions of \S 3.9 and let $\gamma $ be a closed
curve contained in $U$, then
${\hat I}n (0; f\circ \gamma )=\sum_{\partial _f(a)\ne 0}
{\hat I}n (a; \gamma )\partial _f(a)$.}
\par {\bf Proof.} There is the equality ${\hat I}n(0;f\circ \gamma )
=\int_{\zeta \in \gamma }d Ln (f(\zeta ))=$
$\int_0^1d Ln (f\circ \gamma (s))=$ $\int_{\gamma }f^{-1}(\zeta )df(\zeta )$.
Let $\partial _f(a)=n\in \bf N$, then
$$f^{-1}(a)f'(a).S=
\sum_{l,k; n_1+...+n_k=\partial _f(a), 0\le n_j\in {\bf Z}, j=1,...,k}
(z-a)^{n_1} g_{S,l,k,1;n_1,...,n_k}(z)$$
$$(z-a)^{n_2}g_{S,l,k,2;n_1,...,n_k}(z)...
(z-a)^{n_k} g_{S,l,k,k;n_1,...,n_k}(z),$$
where $g_{S,l,p,k;n_1,...,n_k}(z)$ are
quaternion holomorphic functions of $z$ on $U$
such that $g_{S,l,p,k;n_1,...,n_k}(a)\ne 0$, $S\in \{ I,J,K,L \} $,
where $l=1,...,m$, $1\le m\le 4^{\partial _f(a)}$ (see \S \S 2.8,
3.7, 3.22, 3.28), since each term
$\xi (z) (v-v_0)^{n_1}(w-w_0)^{n_2}(x-x_0)^{n_3}(y-y_0)^{n_4}$
with $n_1+...+n_4\ge \partial _f(a)$, $n_j\ge 0$, has such decomposition,
where $\xi (z)$ is a quaternion holomorphic function on a neighbourhood
of $a$ such that $\xi (a)\ne 0$.
Suppose $\psi $ is a closed curve such that ${\hat I}n(p,\psi )=
2\pi nM$, $|M|=1$, $M\in \bf H_i$, $0\ne n\in \bf Z$.
Then we can define a curve $\psi ^{1/n}=:\omega $ as a closed curve
for which ${\hat I}n(p,\omega )=2\pi M$ and $\omega ([0,1])\subset
\psi ([0,1])$. Then we call $\omega ^n=\psi $. That is,
${\hat I}n(p,\psi ^{1/n})={\hat I}(p,\psi )/n$. The latter formula
allows an interpretation also when ${\hat I}n (p,\psi )/n$ is
equal to $2\pi qM$, where $0\ne q\in \bf Q$. That is, a curve
$\psi ^{1/n}$ can be defined for each $0\ne n\in \bf Z$.
This means that $\gamma $ can be presented
as union of curves $\omega _j$ for each of which there exists
$n_j\in \bf N$ such that $\omega _j^{n_j}$ is a closed curve.
Using Theorem 3.9 for each $a\in U$ with $\partial _f(a)\ne 0$, also
using the series given above we can find a finite family of
$\omega _j$ for which one of the terms in the series is not less,
than any other term. We may also use small homotopic
deformations of $\omega _j$ satisfying the conditions of Theorem 3.9
such that in the series one of the terms is greater
than any other for almost all points on $\omega _j$.
Such deformation is permitted, since otherwise two terms would coincide
on an open subset of $U$, that is impossible.
Considering such series, Formulas $(2.6, 2.7)$ and using Theorem 3.27.2
we get the statement of this theorem.
\par {\bf 3.31. Theorem.} {\it If $f$ has an essential
singularity at $a$, then $cl (f(V))={\hat {\bf H}}$
for each $V\subset dom (f)$, $V=U\setminus \{ a \} $,
where $U$ is a neighbourhood of $a$.}
\par {\bf Proof.} Suppose that the statement of this theorem
is false, then there would exist $r>0$ and $m>0$ and a quaternion
$A\in \bf H$ such that $f$ is $z$-analytic in $B(a,0,r,{\bf H})\setminus
\{ a \} $ and $|f(z)-A|\ge m$ for each $z$ such that
$0<|z-a|<r$. If $\infty \notin cl (f(V))$, then there exists
$R>0$ such that $A\notin cl (f(V))$ for each $|A|>R$.
Therefore, the function $[f(z)-A]^{-1}$ is quaternion
holomorphic in $B(a,0,r,{\bf H})\setminus \{ a \} $.
Hence $[f(z)-A]^{-1}=\sum_k(p_k,(z-a)^k)$, where
in this sum $k=(k_1,...,k_{m(k)})$ with $k_j\ge 0$
for each $j=1,...,m(k)\in \bf N$, $p_k$ are finite sequences
of coefficients for $[f(z)-A]^{-1}$ as in \S 3.22.
If $D_z^n([f(z)-A]^{-1})|_{z=a}=0$ for each $n\ge 0$, then
$[f(z)-A]^{-1}=0$ in a neighbourhood of $a$.
Therefore, $[f(z)-A]^{-1}=\sum_{n_1+...+n_l=n}g_1z^{n_1}...
g_lz^{n_l}$ for some $n$ such that $0\le n\in \bf N$,
$n_j\ge 0$ for each $j=1,...,l\in \bf N$,
each $g_j$ is a quaternion holomorphic function (of $z$).
Consequently, taking inverses of both sides
$[f(z)-A]$ and $(\sum_{n_1+...+n_l=n}g_1z^{n_1}...
g_lz^{n_l})^{-1}$ and comparing their expansion
series we see that finite sequences $b_k$ of expansion
coefficients for $f$ have the property $b_k=0$ for each
$\eta (k)<-n$. This contradicts the hypothesis
and proves the theorem.
\par {\bf 3.32. Definition.} Let $a$ and $b$ be two points in $\bf H$
and $\theta $ be a stereographic mapping of the unit four dimensional
real sphere $S^4$ on $\bf \hat H$. Then $\chi (a,b):=
|\phi (a)-\phi (b)|_{\bf R^5}$ is called the chordal metric,
where $\phi :=\theta ^{-1}: {\bf \hat H}\to \bf S^4$,
$S^4$ is embedded in $\bf R^5$ and $|*|_{\bf R^5}$ is the Euclidean
distance in $\bf R^5$.
\par {\bf 3.32.1. Theorem.} {\it Let $U$ be an open
region in $\bf \hat H$, $\{ f_n: n \in {\bf N} \} $ be a sequence
of functions meromorphic on $U$ tending uniformly in $U$
to $f$ relative to the chordal metric. Then either $f$ is the constant
$\infty $ or else $f$ is meromorphic on $U$.}
\par {\bf 3.32.2. Theorem.} {\it Let $\{ f_k: k\in {\bf N} \} $
be a sequence of meromorphic functions on an open subset $U$ in $\bf \hat H$,
which tends uniformly in the sence of the chordal metric in $U$
to $f$, $f\ne const $. If $f(a)=b$ and $r>0$ are such that
$B(a,r,{\bf H})\subset U$ and $f(z)\ne b$ for each $z\in B(a,r,{\bf H})
\setminus \{ a \} $, then there exists $m\in \bf N$ such that
the value of the valence of $f_k|_{B(a,r,{\bf H})}$ at $b$
is $n(b;f)=n(a;f)$ for each $k\ge m$.}
\par {\bf 3.32.3. Note.} The proofs of these theorems are
formally similar to the proofs of VI.4.3 and 4.4 \cite{heins}.
Theorem 3.32.2 is the quaternion
analog of the Hurwitz theorem.
There are also the following quaternion analogs of
the Mittag-Leffler and Weierstrass theorems.
Their proofs are similar to those for Theorems VIII.1.1 and
1.2 respectively. Nevertheless the second part of the Weierstrass
theorem is not true in general because of noncommutativity of $\bf H$,
that is, a function $h\in {\bf M}(U)$
with $\partial _h=\partial $ is not necessarily
representable as $h=fg$, where $g$ is quaternion holomorphic on $U$
and $f$ is another marked function $f\in {\bf M}(U)$ such that
$\partial _f=\partial $.
In the proofs ordered products of more elementary polynomial
functions and in particular linear terms $(z-b_k)$
have to be considered as in \S 3.28, using Theorems
3.17 and 3.22. Theorem 3.33.2 is not true in general
without condition of right superlinearity
(or left superlinearity) of the superdifferential, for example,
the function $f(z)=xK$ serves as a counterexample,
where $z=vI+wJ+xK+yL$, $v$, $w$, $x$ and $y\in \bf R$, $z\in \bf H$.
\par {\bf 3.33. Theorem.} {\it Let $U$ be a nonempty proper open
subset of $\bf \hat H$, let $A\subset U$ not containing any cluster point
in $U$. Let there be a function $g_b\in {\bf M}({\bf \hat H})$
for each $b\in A$ having a pole at $b$ and no other.
Then there exists $f\in {\bf M}(U)$ quaternion holomorphic on
$U\setminus B$ and having the same principal part at $b$ as $g_b$.
If $f$ is such a function, then each other such function
is the function $f+g$, where $g$ is quaternion holomorphic on $U$.}
\par {\bf 3.33.1. Theorem.} {\it Let $U$ be a proper nonempty
open subset of $\bf \hat H$. Let $\partial : U\to \bf Z$ be a
function such that $ \{ \partial (z)\ne 0 \} $ does not have
a cluster point in $U$.
Then there exists $f\in {\bf M}(U)$ such that $\partial _f=\partial $.}
\par {\bf 3.33.2. Theorem.} {\it Let $U$ be an open region
in $\bf H$ and $f$ be a function quaternion holomorphic on $U$
with a right superlinear superdifferential on $U$.
Suppose $f$ is not constant and $B(a,r,{\bf H})\subset U$,
where $0<r<\infty $. Then $f(B(a,r,{\bf H}))$ is a neighbourhood
of $f(a)$ in $\bf H$.}
\par {\bf 3.34. Remarks.} For calculating
expansion coefficients $b_k$ of a function $f$ quaternion
holomorphic on $U\setminus \{ z_0 \} $, where $U$ is open in $\bf H$,
it is possible to use the residues
$res [(f(z)(z-z_0)^l)^{(n)}.(S_{j_1},...,S_{j_n})]$, where
$S_j\in \{ I,J,K,L \} $, $0\le l\in \bf Z$, $0\le n \in \bf Z$.
But the system of equations for each $b_k=(b_{k,1},...,b_{k,m})$
is nonlinear in general. The calculation of a residue of
a term $(b_k,z^k)$ along the closed curve $\gamma (s)=r\exp (2\pi sM)$
(or $\psi $ homotopic to it and satisfying conditions of Theorem 3.9)
with $|M|=1$, $M\in \bf H_i$, $s\in [0,1]$, reduces to a calculation
of a $\bf R$-linear combination of integrals of the form \\
$\int_0^1\exp (2\pi sn_1M_1)...\exp (2\pi s n_lM_l)ds A,$
where $n_1,...,n_l\in \bf Z$, $n_1+...+n_l=0$, $M_j:={\tilde S_j}MS_j,$
$S_j\in \{ I,J,K,L \} $ for each $j=1,...,l$, $A\in \{ J,K,L \} $.
The case of $f\in \mbox{ }_lC^{\omega }(U,{\bf H})$ is trivial due
to Corollary 3.25.
\par For several quaternion
variables a multiple quaternion line integral
${\bf I}:=\int_{\gamma _n}(...(\int_{\gamma _1}f(\mbox{ }^1z,...,\mbox{ }^nz)
d\mbox{ }^1z)...)d\mbox{ }^nz$ may be naturally considered
for rectifiable curves $\gamma _1$,...,$\gamma _n$ in $\bf H$.
If $\gamma _j=r_j\exp (2\pi s_jM_j)$ with $0<r_j<\infty $,
$s_j\in [0,1]$ and $[M_k,M_j]=0$ commute for each $k,j =1,...,n$,
then this integral $\bf I$ does not depend on the order
of integration for $f\in C^0(U,{\bf H})$, where
$U$ is an open subset in $\bf {\hat H}^n$ and $\gamma _j\subset \mbox{ }^jU$
for each $j$, $U=\mbox{ }^1U\times ...\times \mbox{ }^nU$,
$\mbox{ }^jU$ is an open subset in $\bf H$.
Therefore, there is the natural generalization
of Theorem 3.9 for several quaternion variables:
$$(3.9')\quad (2\pi )^nf(z_0)=$$
$$(\int_{\psi _n}(...(
\int_{\psi _1}f(\mbox{ }^1\zeta ,...,\mbox{ }^n\zeta )(\mbox{ }^1\zeta -
\mbox{ }^1z_0)^{-1}d\mbox{ }^1\zeta )M_1^{-1})...)
(\mbox{ }^n\zeta - \mbox{ }^nz_0)^{-1}d\mbox{ }^n\zeta )M_n^{-1} $$
for the corresponding $U=\mbox{ }^1U\times ...\times \mbox{ }^nU,$
where $\psi _j$ and $\mbox{ }^jU$ satisfy conditions of Theorem 3.9
for each $j$ and $f$ is a continuous quaternion holomoprhic function on $U$.

\thanks{
Addresses: Sergey V. L\"udkovsky, Theoretical Department,
Institute of General Physics,
Str. Vavilov 38, Moscow, 119 991 GSP-1, Russia.\\
E-mail: ludkovsk@fpl.gpi.ru; \\
Fred van Oystaeyen, Department of Mathematics and Computer
Sciences, University of Antwerpen, UIA.
Universiteitsplein 1, 2610 Antwerpen, Belgium.
E-mail: francin@uia.ua.ac.be \\
{\underline {Acknowledgment}}. Both authors thank the European
Science Foundation for support through the Noncommutative Geometry
(NOG) project.}
\end{document}